\documentclass[11pt,a4paper,reqno]{amsart}
\usepackage{amssymb,amsmath,amsthm}
\usepackage[utf8]{inputenc}
\usepackage[T1]{fontenc}
\usepackage{enumerate}
\usepackage[all]{xy}
\usepackage{fullpage}
\usepackage{comment}
\usepackage{array}
\usepackage{longtable}
\usepackage{stmaryrd}
\usepackage{mathrsfs} 
\usepackage{xcolor}
\usepackage{mathtools}

\def\gcd{\mathrm{gcd}}
\def\deg{\mathrm{deg}}

\def\dim{\mathrm{dim}}

\usepackage{hyperref}
\hypersetup{
    colorlinks=true,
    linkcolor=blue,
    filecolor=red,  
    citecolor=green,
    urlcolor=cyan,
    pdftitle={Orthogonal Sets},
    pdfpagemode=FullScreen,
    }

\usepackage{cleveref}
\crefformat{section}{\S#2#1#3}
\crefformat{subsection}{\S#2#1#3}
\crefformat{subsubsection}{\S#2#1#3}
\usepackage{enumitem}
\usepackage{tikz}

\newtheorem*{claim*}{Claim}

\newtheorem{theorem}{Theorem}[section]
\newtheorem{lemma}[theorem]{Lemma}

\newtheorem{proposition}[theorem]{Proposition}
\newtheorem{example}[theorem]{Example}

\theoremstyle{definition}
\newtheorem{definition}[theorem]{Definition}

\theoremstyle{remark}
\newtheorem{remark}[theorem]{Remark}
\newtheorem{question}[theorem]{Question}

\makeatletter
\newcommand{\subalign}[1]{%
  \vcenter{%
    \Let@ \restore@math@cr \default@tag
    \baselineskip\fontdimen10 \scriptfont\tw@
    \advance\baselineskip\fontdimen12 \scriptfont\tw@
    \lineskip\thr@@\fontdimen8 \scriptfont\thr@@
    \lineskiplimit\lineskip
    \ialign{\hfil$\m@th\scriptstyle##$&$\m@th\scriptstyle{}##$\hfil\crcr
      #1\crcr
    }%
  }%
}
\makeatother

\numberwithin{equation}{section}

\title{Counting Problems for Orthogonal Sets and Sublattices in Function Fields}
\author{Noy Soffer Aranov}
\email{noy.sofferaranov@tugraz.at}
\address{Graz University of Technology, Institute of Analysis and Number Theory, 8010 Graz, Austria}

\author{Angelot Behajaina}
\email{angelot.behajaina@univ-lille.fr}
\address{Univ. Lille, CNRS, UMR 8524 - Laboratoire Paul Painlevé, F-59000 Lille, France}

\begin{document}

\maketitle
\begin{abstract}
Let $\mathcal{K}=\mathbb{F}_q((x^{-1}))$. Analogous to orthogonality in the Euclidean space $\mathbb{R}^n$, there exists a well-studied notion of ultrametric orthogonality in $\mathcal{K}^n$. In this paper, we extend the work of \cite{AB24} on counting problems related to orthogonality in $\mathcal{K}^n$. For example, we resolve an open question posed in \cite{AB24} by bounding the size of the largest ``orthogonal sets'' in $\mathcal{K}^n$. Furthermore, using similar ideas and techniques, we investigate analogues of Hadamard matrices over $\mathcal{K}$. Finally, we also use ultrametric orthogonality to compute the number of sublattices of $\mathbb{F}_q[x]^n$ with a certain geometric structure, and to determine the number of orthogonal bases of a sublattice in $\mathcal{K}^n$. The resulting formulas depend crucially on successive minima.
\end{abstract}

\section{Introduction}
Erd{\H{o}}s asked the following question: given integers $2 \leq \ell \leq k$, what is the maximum size $\alpha_n^{(k,\ell)}$ of a subset $S \subseteq \mathbb{R}^n \setminus \{{\bf 0}\}$ satisfying the following property: for any $E \subseteq S$ of size $k$, there exists a subset $F \subseteq E$ of size $\ell$ such that vectors in $F$ are pairwise orthogonal\footnote{with respect to the standard scalar product.}? Over the past few decades, this question has attracted considerable interest; however, it remains widely open. Regarding known results, it is clear that $\alpha_n^{(2,2)}=n$. F\"uredi and Stanley \cite{FS92} showed that $\alpha_2^{(k,2)}=2k-2$. Rosenfeld \cite{Ros95} proved that $\alpha_n^{(3,2)}=2n$, which was later reproved by Deaett \cite{Dea11}. See \cite{AS99} for further related results.

Ahmadi and Mohammadian \cite{AM} investigated an analogue of Erd{\H{o}}s' question in the finite field setting -- specifically, in $\mathbb{F}_q^n$ equipped with a symmetric non-degenerate bilinear form. For example, they computed the corresponding $\alpha_{n}^{(2,2)}$ and gave an upper bound on $\alpha_{n}^{(3,2)}$. The later was improved by Mohammadian and Petridis \cite{MP22}.

In \cite{AB24}, we investigated variants of Erd{\H{o}}s' problem   over a discrete valued field. In this setting, the notion of ‘‘orthogonality'' is strongly inspired (see Definition \ref{def:orthvect} below) by the following characterization in the classical setting:
$$
{\bf v},{\bf w} \in \mathbb{R}^n\textrm{ are orthogonal if and only if }\Vert \lambda {\bf v}+ \mu {\bf w}\Vert_2^2=\Vert \lambda {\bf v} \Vert_2^2+ \Vert \mu {\bf w} \Vert_2^2,\,\textrm{ for all }\lambda,\mu \in \mathbb{R},
$$ where $\Vert\cdot \Vert_2$ denotes the $\ell_2$-norm. Now fix a discrete valued field $\mathcal{K}$ with a {\bf finite residue} field $\kappa=\mathbb{F}_q$ (where $q$ is a prime power), a valuation $\nu: \mathcal{K} \rightarrow \mathbb{Z}$ and a maximal ideal\footnote{that is $\mathfrak{m}=\{ \lambda \in \mathcal{K} \mid \nu(\lambda)>0\}$.} $\mathfrak{m}$. The absolute value $|\cdot|$ on $\mathcal{K}$ is given by:
$$
|\lambda|=q^{-\nu(\lambda)},\,\,\,\textrm{for all}\,\, \lambda \in \mathcal{K}.
$$
This extends to the \emph{infinity norm} on $\mathcal{K}^n$ ($n \geq 1$) by:
$$
\Vert{\bf v}\Vert=\max_{1 \leq i \leq n}\vert v_i\vert ,\,\, \textrm{for all}\,\, {\bf v}=(v_1,\dots,v_n) \in \mathcal{K}^n.
$$ 
The norm and absolute value on $\mathcal{K}^n$ satisfy the ultrametric inequality.

\begin{lemma}[Ultrametric Inequality]
\label{lem:UMIneq}
    Let $n\in \mathbb{N}$ and let $\mathcal{K}$ be a discrete valued field with finite residue. Then, for every $\mathbf{u}_1,\dots,\mathbf{u}_{\ell}\in \mathcal{K}^n$, we have 
    \begin{equation}
    \label{eqn:UMIneq}
        \Vert \mathbf{u}_1+\dots+\mathbf{u}_{\ell}\Vert\leq \max_{i=1,\dots,\ell}\Vert \mathbf{u}_i\Vert.
    \end{equation}
    Moreover, if there exists $i_0\in \{1,\dots,\ell\}$, such that for every $i\in \{1,\dots,\ell\}\setminus\{i_0\}$, we have $\Vert \mathbf{u}_i\Vert<\Vert \mathbf{u}_{i_0}\Vert$, then, the inequality in \eqref{eqn:UMIneq} is an equality.
\end{lemma}

\begin{example}\label{ex:dvr}
Here are some natural examples of discrete valued fields with finite residue: \newline
\noindent $\bullet$ Let $\mathcal{K}=\mathbb{F}_q((x^{-1}))$ be the field of Laurent series in $x^{-1}$, and let $\nu$ be the valuation defined by $\nu(\sum_{i}a_ix^{-i})=\min(i \mid a_i\neq 0)$. In this case, $\mathfrak{m}=x^{-1} \mathbb{F}_q[[x^{-1}]]$ and $\kappa=\mathbb{F}_q$. 
\vskip 1mm
\noindent $\bullet$ Let $\mathcal{K}=\mathbb{Q}_p$ ($p$ a prime number) be the field of $p$-adic numbers, and let $\nu$ be the $p$-adic valuation. In this case, $\mathfrak{m}=p \mathbb{Z}_p$ and $\kappa=\mathbb{F}_p$.
\end{example}
\noindent

\begin{definition}\label{def:orthvect} 
Let $n\in \mathbb{N}$.\newline
\noindent $\bullet$
We say that vectors $\mathbf{u}_1,\dots,\mathbf{u}_\ell \in \mathcal{K}^n \setminus \{{\bf 0}\}$ ($\ell \geq 1$) are \emph{orthogonal} if 
    $$\Vert \lambda_1 \mathbf{u}_1+\dots + \lambda_\ell \mathbf{u}_\ell\Vert=\max_{1 \leq i \leq \ell}\Vert \lambda_i {\bf u}_i \Vert,$$
    for all $\lambda_1,\dots,\lambda_\ell \in \mathcal{K}$. 
\vskip 1mm
\noindent $\bullet$ We say that vectors $\mathbf{u}_1,\dots,\mathbf{u}_\ell \in \mathcal{K}^n \setminus \{{\bf 0}\}$ ($\ell \geq 1$) are \emph{weakly orthogonal} if they are pairwise orthogonal. 
\end{definition}

\begin{remark}
    Weakly orthogonal vectors need not be orthogonal. For example, consider $$\mathbf{u}_1=\begin{pmatrix}
        1\\
        0\\
        0
    \end{pmatrix},\mathbf{u}_2=\begin{pmatrix}
        0\\
        1\\
        0
    \end{pmatrix},
        \mathbf{u}_3=\begin{pmatrix}
            1\\
            1\\
            0
        \end{pmatrix}\in \mathbb{F}_q((x^{-1}))^3.$$ Indeed, every pair of $\mathbf{u}_i,\mathbf{u}_j$ ($i\neq j$) is orthogonal. However, since $\mathbf{u}_1+\mathbf{u}_2-\mathbf{u}_3=0$, the vectors $\mathbf{u}_1,\mathbf{u}_2,\mathbf{u}_3$ are not orthogonal.
\end{remark}
Recall the following variants of $\alpha_n^{(k,\ell)}$ of the Erd{\H{o}}s' problem that we are investigating. They were introduced in \cite{AB24}. The \emph{unit sphere} in $\mathcal{K}^n$ is defined by $\mathbb{B}_n=\{{\bf v}\in \mathcal{K}^n \mid \Vert{\bf v}\Vert=1\}$, so that $\mathbb{B}_n/\mathfrak{m}^n\cong \mathbb{F}_q^n$.

\begin{definition}\label{ref:strongorthsetdef}
Let $2 \leq \ell \leq k$. 
\vskip 1mm
\noindent $\bullet$ A subset $S\subseteq \mathcal{K}^n\setminus \{{\bf 0}\}$ is \emph{$(k,\ell)$-weakly orthogonal} if any $E\subseteq S$ of size $k$ contains a weakly orthogonal subset $F \subseteq E$ of size $\ell$. We define:
    \begin{align*}
    \Delta_{n}^{(k,\ell)}&=\max\{\vert S\vert  \mid S \subseteq \mathbb{B}_n\,\,\textrm{is}\,\,(k,\ell)\textrm{-weakly orthogonal}\}.
    \end{align*} 
\vskip 1mm
\noindent $\bullet$ A subset $S\subseteq \mathcal{K}^n\setminus \{{\bf 0}\}$ is \emph{$(k,\ell)$-orthogonal} if any $E\subseteq S$ of size $k$ contains an orthogonal subset $F \subseteq E$ of size $\ell$. We define:
    \begin{align*}
    \Theta_{n}^{(k,\ell)}&=\max\{\vert S\vert  \mid S \subseteq \mathbb{B}_n\,\,\textrm{is}\,\,(k,\ell)\textrm{-orthogonal}\}.
    \end{align*}
\end{definition}
In \cite{AB24}, we provided asymptotic bounds for $\Delta_{n}^{(k,\ell)}$ and compared it with $\Theta_n^{(k,\ell)}$. In \cite{AB24}, we also obtained partial asymptotic results for $\Theta_n^{(k,\ell)}$. Here, we strengthen these results and provide a complete picture. To do so, we also need to study the maximum size $
{\rm Ind}_q(n,k,\ell)$ of a subset $S \subset \mathbb{F}_q^n \setminus \{{\bf 0} \}$ satisfying the following property: any subset $X \subset S$ of size $k$ satisfies $\dim(\langle X \rangle) \geq \ell$; such a $S$ is called \emph{$(k,\ell)$-independent}. The investigation of independent sets has appeared in earlier works, such as \cite{Tal61,Hir83,DMMS,DMM,TV09}.

\begin{theorem}[Maximum size of orthogonal sets]\label{thm:orthogset} Let $2 \leq \ell \leq k$. 
\begin{enumerate}
\item (\cite[Question 3.28]{AB24}). We have 
\begin{equation}\label{eq:boundalphakln}
            \left\lfloor \frac{(k-1)(q-1)}{q^{\ell-1}-1}\right\rfloor\frac{q^n-1}{q-1}\leq \Theta_{n}^{(k,\ell)}\leq \left\lfloor \frac{k-1}{q^{\ell-1}-1} (q^n-1)\right\rfloor.
        \end{equation} 
         As a consequence, if $\frac{q^{\ell-1}-1}{q-1}$ divides $k-1$, then $\Theta_{n}^{(k,\ell)}=\frac{k-1}{q^{\ell-1}-1}\times(q^n-1)$. Moreover, $$
\lim_{k \rightarrow \infty }\frac{\Theta_n^{(k,\ell)}}{k-1}=\frac{q^n-1}{q^{\ell-1}-1}. 
$$ \label{thm:orthogsetsizemaxTh} 
\item (\cite[Question 3.17]{AB24}). $\Theta_n^{(k,\ell)}=\Delta_n^{(k,\ell)}$ if and only if $\ell=2$.\label{thm:orthogsetthetadeltcomp}
\item Suppose that $\Theta_{n}^{(k,\ell)}={\rm Ind}_q(n,k,\ell)$. Then $k \leq q^{\ell-1}$. \label{thm:orthogsettheind}
\end{enumerate}       
\end{theorem}

\begin{remark}
In \cite[Theorem 1.17 (2)]{AB24}, the statement ``Moreover, if the equality holds, then $k \leq q^{\ell-2}+1$'', which is \cite[Proposition 3.24]{AB24}, is false. In fact, the mistake appears in the proof of Proposition 3.24. Indeed, the conclusion ``Therefore, by Proposition 3.14., $\gamma_n({\bf v})$ belongs to every $\ell-1$-dimensional {\bf subset} of $\gamma_n(X)$'' should be ``... $\gamma_n({\bf v})$ belongs to every $\ell-1$-dimensional {\bf subspace} of $\langle \gamma_n(X) \rangle$. Theorem \ref{thm:orthogset}(\ref{thm:orthogsettheind}) gives the correct bound. 
\end{remark}
Under some natural assumptions, we can characterize maximum (resp. weakly) orthogonal sets. To do so, we define $\gamma_n:\mathbb{B}_n \rightarrow \mathbb{F}_q^n\setminus \{{\bf 0} \}$ and $\rho_n:\mathbb{F}_q^n\setminus \{{\bf 0} \} \rightarrow \mathbb{P}^{n-1}(\mathbb{F}_q)$ as the reduction modulo $\mathfrak{m}$ and the canonical projection, respectively.

\begin{theorem}[Maximum orthogonal sets]\label{thm:structumaxsize}
Let $2 \leq \ell \leq k$.
\begin{enumerate}
\item Suppose that $\ell-1$ divides $k-1$. Then the maximum $(k,\ell)$-weakly orthogonal sets are exactly the sets $\mathcal{F} \subset \mathbb{B}_n$ such that
$$
\left|\mathcal{F} \cap (\rho_n \circ \gamma_n)^{-1}({\bf v})\right|=\frac{k-1}{\ell-1},
$$ for all ${\bf v} \in \mathbb{P}^{n-1}(\mathbb{F}_q)$.\label{thm:structumaxsizeWOcase}
\item Suppose that $\frac{q^{\ell-1}-1}{q-1}$ divides $k-1$. Then the maximum $(k,\ell)$-orthogonal sets are exactly the sets $\mathcal{F} \subset \mathbb{B}_n$ such that 
$$
\left|\mathcal{F} \cap (\rho_n \circ \gamma_{n})^{-1}({\bf v})\right|=\frac{(k-1)(q-1)}{q^{\ell-1}-1},
$$ for all ${\bf v} \in \mathbb{P}^{n-1}(\mathbb{F}_q)$.\label{thm:structumaxsizeOrcase}
\end{enumerate}
\end{theorem}

In order to prove the above results, we used the following tool which connects between orthogonal sets in $\mathcal{K}^n$ and independent sets in $\mathbb{F}_q^n$.
\begin{proposition}{\cite[Proposition 3.14]{AB24}}
\label{prop:OrtIndConn}
    A set $S\subseteq \mathbb{B}_n$ is orthogonal if and only if the set $\gamma_n(S)\subseteq \mathbb{F}_q^n$ is linearly independent. Moreover, $S\subseteq \mathbb{B}_n$ is $(k,\ell)$-orthogonal if and only if for every $X\subseteq S$ with $\vert X\vert=k$, we have $\dim(\operatorname{span}_{\mathbb{F}_q}(\gamma_n(S))\geq \ell$.
\end{proposition}
Applications of Proposition \ref{prop:OrtIndConn} extend far beyond analogues of the Erd{\H{o}}s' problem. To illustrate this, we study Hadamard matrices in the function field setting. Moreover, we count sublattices of $\mathbb{F}_q[x]^n$ of a given height. 

The second part of our work focuses then on an application of Proposition \ref{prop:OrtIndConn} to study an analogue of Hadamard matrices in polynomial rings over finite fields. Recall that Hadamard matrices over $\mathbb{R}$ are $n\times n$-matrices ($n \geq 1$) whose entries belong to $\{-1,+1\}$ and whose columns are pairwise orthogonal. These matrices have numerous applications, particularly in error correcting codes, due to the fact that their determinant is maximal among all matrices with entries in $\{-1,+1\}$. Now let $\mathcal{K}=\mathbb{F}_q((x^{-1}))$ be as in Example \ref{ex:dvr}. In order to define an analogue of Hadamard matrices over $\mathcal{K}$, let $\mathcal{R}=\mathbb{F}_q[x]$ be the set of all polynomials over $\mathbb{F}_q$. Moreover, for $\ell\geq 0$, let $\mathcal{R}_{\leq \ell}$ denote the set of polynomials of degree at most $\ell$. 

\begin{definition}
Let $n \geq 1$, $\ell \geq 0$ and $m \leq \ell n$.
   \begin{enumerate} 
    \item A matrix $A \in M_{n\times n}(\mathcal{R}_{\leq \ell})$ is called \emph{$(\ell,q^m)$-Hadamard} if the columns of $H$ form an orthogonal set and $\vert\det(H)\vert=q^m$. The set of $(\ell,q^m)$-Hadamard matrices is denoted by $\mathcal{H}_n(\ell,q^m)$.
    \item A matrix $H\in M_{n\times n}(\mathcal{R}_{\leq \ell})$ is called \emph{$(\ell,q^m)$-weakly Hadamard} if the columns of $H$ form a weakly orthogonal set and $\vert \det(H)\vert=q^m$. The set of $(\ell,q^m)$-weakly Hadamard matrices is denoted by $\mathcal{W}_n(\ell,q^m)$.
   \end{enumerate}
\end{definition}
\noindent
We determine the number of $(\ell,q^m)$-Hadamard (respectively weakly Hadamard) matrices and establish an analogue of a result by Sylvester \cite{Syl}, which concerns the stability of the Hadamard property under Kronecker product. 

\begin{theorem}\label{thm:hadamard} 
Let $n\geq 1$, $\ell \geq 0$ and $m \leq \ell n$. Then 
\begin{enumerate}
\item\label{thm:hadamardcounting} $
|\mathcal{H}_n(\ell,q^m)|=q^m \prod_{i=0}^{n-1}(q^n-q^i)\left\vert\left\lbrace(\ell_1,\dots,\ell_n) \in \{0,\dots,\ell\}^n \Bigm\vert \sum_{i=1}^n \ell_i=m\right\rbrace \right\vert,
$
and
\item \label{thm:weakHadamarkCounting}$
|\mathcal{W}_n(\ell,q^m)|=q^m (q-1)^n n! \binom{\frac{q^n-1}{q-1}}{n} \left\vert\left\lbrace(\ell_1,\dots,\ell_n) \in \{0,\dots,\ell\}^n \Bigm\vert \sum_{i=1}^n \ell_i=m\right\rbrace\right\vert.
$
\end{enumerate}
\end{theorem}

\begin{theorem}
    Let $n_1,n_2 \geq 1$, $\ell_1,\ell_2 \geq 0$ and $m_1 \leq \ell_1 n_1,m_2 \leq \ell_2 n_2$. Given an $(\ell_1,q^{m_1})$-Hadamard (respectively weakly Hadamard) matrix $A\in M_{n_1\times n_1}(\mathcal{R}_{\leq \ell_1})$ and an $(\ell_2,q^{m_2})$-Hadamard (respectively weakly Hadamard) matrix $B\in M_{n_2\times n_2}(\mathcal{R}_{\leq \ell_2})$, the Kronecker product $A\otimes B$ is an $(\ell_1+\ell_2,q^{m_1n_2+m_2n_1})$-Hadamard  (respectively weakly Hadamard) matrix.
\label{thm:hadamardKroneckerprod}
\end{theorem}
The final part of our work uses Proposition \ref{prop:OrtIndConn} to examine the structure of sublattices in $\mathcal{K}^n$ ($n \geq 1$) where $\mathcal{K}=\mathbb{F}_q((x^{-1}))$ is  as in Example \ref{ex:dvr}. Specifically, we compute the number of sublattices of $\mathbb{F}_q[x]^n$, and determine the number of orthogonal bases of a sublattice in $\mathcal{K}^n$. Orthogonality, as shown by the results in \cite{Mah,RW,BK,Ara}, plays a crucial role in our investigation.

\begin{remark}
    Here, we focus on $\mathcal{K}=\mathbb{F}_q((x^{-1}))$, as sublattices in $\mathcal{K}^n$ admit a ``Minkowski basis'' \cite{Mah,Ara}. This property does not hold for general discrete valued fields and for general global function fields. 
\end{remark}

\begin{definition}
Let $n\in\mathbb{N}$, and let $1\leq k\leq n$.
\vskip 1mm
\noindent $\bullet$ A \emph{lattice} in $\mathcal{K}^n$ is a subgroup of the form $\Lambda= g \times \mathcal{R}^n$ with $g \in {\rm GL}_n(\mathcal{K}^n)$.
\vskip 1mm
\noindent $\bullet$ A subgroup $\Lambda \leq \mathcal{K}^n$ is called a \emph{$k$-sublattice} or a \emph{sublattice of dimension $k$} if $\dim(\operatorname{span}_{\mathcal{K}}(\Lambda))=k$ and $\Lambda$ is a lattice in $\operatorname{span}_{\mathcal{K}}(\Lambda)$. Equivalently, there exist $\mathcal{K}$-independent vectors $\mathbf{v}_1,\dots,\mathbf{v}_k\in \mathcal{K}^n$, such that: 
$$\Lambda=\mathcal{R}\mathbf{v}_1\oplus\dots\oplus\mathcal{R}\mathbf{v}_k.$$
In this case, we define $\vert\det(\Lambda)\vert=\Vert \mathbf{v}_1\wedge\dots\wedge\mathbf{v}_k\Vert$.

We denote the space of $k$-sublattices with entries in $\mathcal{K}$ by $\mathcal{L}_{n,k}(\mathcal{K})$. We let $\mathcal{L}_{n,k}(\mathcal{R})$ denote the set of $k$-sublattices $\Lambda=\mathcal{R}\mathbf{v}_1\oplus\dots\oplus\mathcal{R}\mathbf{v}_k$, such that $\mathbf{v}_1,\dots,\mathbf{v}_k\in \mathcal{R}^n$.
\end{definition}
By \cite[Lemma 1.18 and Theorem 1.21]{Ara}, every lattice in $\mathcal{K}^n$ has an orthogonal $\mathcal{R}$-basis; moreover, all such bases correspond to bases of successive minima. These results can be extended to sublattices. To do so, we first provide the following.

\begin{definition}
    Let $\Lambda\leq \mathcal{K}^n$ be a $k$-sublattice ($1 \leq k \leq n$). For $1 \leq i \leq k$, the \emph{$i$-th successive minima} of $\Lambda$ is defined by 
    \begin{align*}
    \lambda_i(\Lambda)&= \inf\left\{ r>0 \mid  \textrm{there exist}\,\,\,\, i\,\,\,\textrm{linearly independent vectors in}\,\,\Lambda\,\,\textrm{of norm}\,\,\leq r \right\} \\
    &=\inf\left\{r>0 \mid \dim\operatorname{span}(\Lambda\cap \{\mathbf{v}:\Vert \mathbf{v}\Vert\leq r\}) \geq i\right\}.
    \end{align*}
By abuse of terminology, an independent set of vectors $\{\mathbf{v}_1,\dots,\mathbf{v}_k\}\subseteq \Lambda$ satisfying $\Vert \mathbf{v}_i\Vert=\lambda_i(\Lambda)$ is called a set of successive minima.
\end{definition}

\begin{proposition}
   \label{thm:AraBase}
    Let $\Lambda \leq \mathcal{K}^n$ be a $k$-sublattice. Then, $\Lambda$ has an orthogonal $\mathcal{R}$-basis consisting of successive minima. Moreover, an $\mathcal{R}$-basis of $\Lambda$ is orthogonal if and only if it is a basis of successive minima. 
\end{proposition}
Motivated by Proposition \ref{thm:AraBase}, we investigate the following question.

\begin{question}
\label{ques:numBases}
How many orthogonal $\mathcal{R}$-bases does a $k$-sublattice $\Lambda \leq \mathcal{K}^n$ possess?
\end{question}
We answer this question, noting that it depends on the multiplicity of the successive minima. Hence, we introduce the following.

\begin{definition}\label{def:ussublatt}
    Let $\boldsymbol{\mu}=(\mu_1,\dots,\mu_m)\in (q^\mathbb{Z})^m$ ($m \geq 1$) satisfy $\mu_1<\dots<\mu_m$, and let $\boldsymbol{s}=(s_1,\dots,s_m)\in \mathbb{N}^m$. A $k$-sublattice $\Lambda \leq \mathcal{K}^n$ is called a \emph{$(\boldsymbol{\mu},\boldsymbol{s})$-sublattice} if: 
    \begin{enumerate}
        \item $\{\lambda_1(\Lambda),\dots,\lambda_k(\Lambda)\}=\{\mu_1,\dots,\mu_m\}$ and,
        \item for each $1 \leq j \leq m$, we have $\vert \{1 \leq i \leq k \mid \lambda_i(\Lambda)=\mu_j\}\vert=s_j$.
    \end{enumerate} 
    We denote the set of $(\boldsymbol{\mu},\boldsymbol{s})$-sublattices of dimension $k$ by $G_{n,k,q}(\boldsymbol{\mu},\boldsymbol{s})$. Furthermore, let $G_{n,k,q}(\boldsymbol{\mu},\boldsymbol{s},\mathcal{R})$ denote the subset of $(\boldsymbol{\mu},\boldsymbol{s})$-sublattices of dimension $k$ contained in $\mathcal{R}^n$.
\end{definition}
Proposition \ref{thm:AraBase} motivates the following definition.
\begin{definition}
An orthogonal basis tuple $(\mathbf{v}_1,\dots,\mathbf{v}_k)$ of a $k$-sublattice $\Lambda$ is said to be \emph{ordered} if $\Vert \mathbf{v}_1\Vert\leq \dots\leq \Vert \mathbf{v}_k\Vert$.
\end{definition}

\begin{theorem}
\label{thm:numOrtbases}
    Let $1\leq k\leq n$, let $\boldsymbol{\mu}\in (q^{\mathbb{Z}})^m$, let $\mathbf{s}\in \mathbb{N}^m$, and let $\Lambda$ be a $(\boldsymbol{\mu},\boldsymbol{s})$-sublattice. Then, the number of orthogonal ordered bases of $\Lambda$ is:
    \begin{equation}
        \begin{cases}
            (q^k-1)\cdots(q^k-q^{k-1})&\textrm{if}\,\,m=1,\\
            q^{\sum_{\ell=1}^ms_{\ell}\sum_{j=1}^{\ell-1}s_j}\prod_{\ell=1}^m\left(\mu_{\ell}^{s_{\ell}\left(\sum_{i=1}^{\ell-1}s_i-\sum_{i=\ell+1}^ms_i\right)}\prod_{i=0}^{s_{\ell}-1}(q^{s_{\ell}}-q^i)\right)&\textrm{otherwise}.
        \end{cases}
    \end{equation}
\end{theorem}
We then apply Theorem \ref{thm:numOrtbases} to compute $\vert G_{n,k,q}(\boldsymbol{\mu},\boldsymbol{s},\mathcal{R})\vert$. This problem can be viewed as an analogue of the Gauss circle problem \cite{Har} in $\mathcal{L}_{n,k}(\mathcal{K})$, which concerns the asymptotic number of lattice points within a ball of large radius, since 
$$\left|\left\{\Lambda\in \mathcal{L}_{n,k}(\mathcal{R}):\vert\det(\Lambda)\vert=q^r\right\}\right|=\sum_{(\boldsymbol{\mu},\boldsymbol{s}):\prod\mu_i^{s_i}=r}\left|G_{n,k,q}(\boldsymbol{\mu},\boldsymbol{s})\right|.$$ 
The Gauss circle problem is a very well studied problem, which has been studied in various settings and dimensions \cite{Har,HarBook,Guy,BKZ,Low,Gath}. Analogues of this question have been studied in finite extensions of function fields \cite{TW} and global fields \cite{HPLat,HPsub}. 

In the function field setting Bagshaw and Kerr \cite[Lemma 6.2] {BK} computed the exact value of $\vert \Lambda\cap B(0,R)\vert$, for every lattice $\Lambda\in \mathcal{L}_n$ and for every $R>0$. We generalize \cite[Lemma 6.2]{BK} to $k$-sublattices. Due to the nature of $\mathcal{K}$ \cite[Theorems 1.18 and 1.21]{Ara} (see also \cite{Mah,Ros}), we obtain exact bounds, whereas in other settings, the bounds are only asymptotic. 

\begin{theorem}
\label{thm:NumG_n,k,q(mu,s)}
Let $1\leq k\leq n$. Then, for every $1\leq m\leq k$, for every $\boldsymbol{\mu}\in (q^{\mathbb{Z}})^m$ with $1 \leq\mu_1<\dots<\mu_m$, and for every $\boldsymbol{s}\in \mathbb{N}^m$, we have:
    \begin{equation}
      \vert G_{n,k,q}(\boldsymbol{\mu},\boldsymbol{s},\mathcal{R})\vert=\begin{cases}
        \mu_1^{kn}&\textrm{if}\,\,m=1,\\
        \frac{\prod_{i=0}^{k-1}(q^n-q^i)}{\prod_{\ell=1}^m\prod_{i=0}^{s_{\ell}-1}(q^{s_{\ell}}-q^i)}\frac{\prod_{\ell=1}^m\mu_{\ell}^{s_\ell\left(n-\sum_{i=1}^{\ell-1}s_i+\sum_{i=\ell+1}^ms_i\right)}}{q^{\sum_{1 \leq j<\ell \leq m}s_js_{\ell}}}&\text{otherwise}.
    \end{cases}
    \end{equation}
\end{theorem}
A refinement of the Gauss circle problem is the primitive Gauss circle problem \cite{Wu}, which concerns the number of primitive lattice points within a ball of large radius. 

\begin{definition}
    A $k$-sublattice $\Lambda \leq \mathcal{R}^n$ is \emph{primitive} if it cannot be expressed as $\Lambda=c\Delta$, where $\Delta<\mathcal{R}^n$ is a $k$-sublattice and $c\in \mathcal{R}$ with $\deg(c)\geq 1$. The set of primitive $k$ dimensional $(\boldsymbol{\mu},\boldsymbol{s})$-sublattices of $\mathcal{R}^n$ is denoted by $\widehat{G}_{n,k,q}(\boldsymbol{\mu},\boldsymbol{s},\mathcal{R})$.
\end{definition}
In \cite{HPsub}, Horesh and Paulin studied the asymptotics of the number of primitive $k$-sublattices with determinant at most $R$ in a global field. Their results pertain to highly general fields, and moreover, address the equidistribution of such sublattices. In this paper, we focus on the function field setting (that is, on $\mathbb{F}_q(x)$) setting, and compute $\vert \widehat{G}_{n,k,q}\left(\boldsymbol{\mu},\boldsymbol{s},\mathcal{R}\right)\vert$. While the results of \cite{HPsub} are asymptotic, we provide an exact value for $\left|\widehat{G}_{n,k,q}\left(\boldsymbol{\mu},\boldsymbol{s},\mathcal{R}\right)\right|$.

\begin{theorem}
\label{thm:PrimLattCount}
    Let $1\leq k\leq n$, let $1\leq m\leq k$, let $\boldsymbol{\mu}\in (q^{\mathbb{Z}})^m$ satisfy $\mu_1<\dots<\mu_m$, and let $\boldsymbol{s}\in \mathbb{N}^m$.  
    \begin{enumerate}
        \item If $\mu_1=1$, then, $\widehat{G}_{n,k,q}(\boldsymbol{\mu},\boldsymbol{s},\mathcal{R})=G_{n,k,q}(\boldsymbol{\mu},\boldsymbol{s},\mathcal{R})$.
        \item If $m=1$ and $\mu_1\neq 1$, then $\vert \widehat{G}_{n,k,q}(\boldsymbol{\mu},\boldsymbol{s},\mathcal{R})\vert=\mu_1^{n_k}\left(1-q^{-(nk-1)}\right)$.
        \item If $m \geq 2$ and $\mu_1 \neq 1$, then:
        $$\vert \widehat{G}_{n,k,q}(\boldsymbol{\mu},\boldsymbol{s},\mathcal{R})\vert=\frac{\prod_{i=0}^{k-1}(q^k-q^i)\times \prod_{j=1}^m\left(\mu_j^{s_j\left(n-\sum_{i=1}^{\ell-1}s_i+\sum_{i=\ell+1}^ms_i\right)}\times \left(1-q^{1-s_j\left(n-\sum_{i=1}^{\ell-1}s_i+\sum_{i=\ell+1}^ms_i\right)}\right)\right)}{q^{\sum_{\ell<j}s_js_{\ell}}\prod_{\ell=1}^m\prod_{i=0}^{s_{\ell}-1}(q^{s_{\ell}}-q^i)}.$$
    \end{enumerate}
\end{theorem}

\subsection*{Organization of the paper} In \cref{sec:prelim}, we establish the preliminary results required for the subsequent sections. In \cref{sec:maxorthset}, we focus on the study of maximum orthogonal sets, providing proofs of Theorem \ref{thm:orthogset} and Theorem \ref{thm:structumaxsize}. In \cref{sec:Hadamard}, we proceed to prove Theorem \ref{thm:hadamard}. Finally, in \cref{sec:SubLattice}, we conclude by proving Theorem \ref{thm:numOrtbases}, Theorem \ref{thm:NumG_n,k,q(mu,s)}, and Theorem \ref{thm:PrimLattCount}.  

\subsection*{Acknowledgements} 
The first author would like to thank Valerii Sopin for introducing her to Hadamard matrices, and thus, motivating sections \cref{sec:Hadamard} and \cref{sec:SubLattice}. The second author is supported by the Labex CEMPI (ANR-11-LABX-0007-01). Both authors would like to thank the anonymous referee, whose comments helped improve the quality of the paper. 

\section{Preliminaries}\label{sec:prelim}
This section, which is primarily based on \cite{AB24} and \cite{Ara}, gathers the basic concepts of orthogonal sets and the geometry of numbers in function fields that we will need later. For that, let $\mathcal{K}$ be a discrete valued field with a finite residue field $\kappa=\mathbb{F}_q$, a valuation $\nu$ and a maximal ideal $\mathfrak{m}$. For $n \geq 1$, recall that $\mathbb{B}_n=\{{\bf v} \in \mathcal{K}^n \mid \Vert {\bf v} \Vert =1 \}$ is the unit sphere. Furthermore, let $\gamma_n:\mathbb{B}_n \rightarrow \mathbb{F}_q^n\setminus \{{\bf 0} \}$ and $\rho_n:\mathbb{F}_q^n\setminus \{{\bf 0} \} \rightarrow \mathbb{P}^{n-1}(\mathbb{F}_q)$ denote the reduction modulo $\mathfrak{m}$ and the canonical projection, respectively. For a subset $S \subset \mathbb{F}_q^n$, the projective subspace $\langle \rho_n(S) \rangle$ of $\mathbb{P}^{n-1}(\mathbb{F}_q)$ is defined as $\rho_n(\langle S \rangle \setminus \{0\})$.

\begin{remark}\label{rmk:subspee}
For $S \subset \mathbb{F}_q^n$, we have $\dim(\langle \rho_n(S) \rangle)=\dim(\langle S \rangle)-1$.
\end{remark}

\subsection{Some lemmas on orthogonality}
The following result is a reformulation of \cite[Proposition 3.14.]{AB24} and relates orthogonality to the reduction map $\rho_n \circ \gamma_n$. 

\begin{lemma}\label{lem:charoforthusingproj}
A set $S=\{{\bf v}_1,\dots,{\bf v}_\ell\} \subset \mathbb{B}_n$ is orthogonal if and only if the vectors $\gamma_n({\bf v}_1),\dots,\gamma_n({\bf v}_\ell)$  are linearly independent. Moreover $S$ is $(k,\ell)$-orthogonal
if and only if for every $S' \subset S$ of size $k$, we have $\dim(\langle (\rho_n \circ \gamma_n)(S')\rangle) \geq \ell-1$.   
\end{lemma}

\begin{proof}
The first statement appears verbatim in \cite[Proposition 3.14.]{AB24}. The second statement follows from \cite[Proposition 3.14.]{AB24} together with Remark \ref{rmk:subspee}.
\end{proof}

For any $S \subseteq \mathbb{B}_n$ and $ {\bf w} \in \mathbb{P}^{n-1}(\mathbb{F}_q)$, define: 
    $$
    t_{{\bf w},S}=\vert \{{\bf v} \in S \mid (\rho_n \circ \gamma_n)({\bf v})={\bf w}\}\vert .
    $$

\begin{lemma}\label{lem:rt}
    Let $2 \leq \ell \leq k$. Then the following are equivalent:
    \begin{enumerate}
    \item $S$ is $(k,\ell)$-orthogonal;
    \item For any subspace $\mathcal{I} \leq \mathbb{P}^{n-1}(\mathbb{F}_q)$ of dimension at most $\ell-2$, we have $\sum_{{\bf w} \in \mathcal{I}}t_{{\bf w},S} \leq k-1$.
    \end{enumerate}\label{lem:ret451}
    \end{lemma}
    
\begin{proof}
By Lemma \ref{lem:charoforthusingproj}, $S$ is $(k,\ell)$-orthogonal, if and only if
\begin{equation}\label{eq:condet1}
\dim(\langle (\rho_n\circ \gamma_n)(S') \rangle)\geq \ell-1,\quad \textrm{for all } S'\subseteq S \textrm{ of size } k.
\end{equation}
This condition is equivalent to the following: 
\vspace{1pt}
\begin{equation}\label{eq:condet2}
\textrm{For any subspace } \mathcal{I} \leq \mathbb{P}^{n-1}(\mathbb{F}_q) \textrm{ with } \dim( \mathcal{I} ) \leq \ell-2, \textrm{ we have } \vert  (\rho_n \circ \gamma_n)^{-1}(\mathcal{I}) \cap S\vert  \leq k-1.
\end{equation}
\vspace{1pt}
\noindent Indeed, assume that \eqref{eq:condet1} holds, and let $\mathcal{I} \leq \mathbb{P}^{n-1}(\mathbb{F}_q)$ be a subspace with $\dim(\mathcal{I})\leq \ell-2$. Define $S'=(\rho_n \circ \gamma_n)^{-1}(\mathcal{I}) \cap S$. Then $\dim(\langle (\rho_n \circ \gamma_n)(S') \rangle)\leq \dim(\mathcal{I}) \leq \ell-2$, and hence \eqref{eq:condet1} implies that $S' \leq k-1$. Therefore, \eqref{eq:condet2} holds. Conversely, assume that \eqref{eq:condet2} holds, and let $S' \subset S$ be a subset of size $k$. Set $\mathcal{I}=\langle (\rho_n \circ \gamma_n)(S') \rangle$. Since $k \leq |S'| \leq |(\rho_n \circ \gamma_n)^{-1}(\mathcal{I})|$, condition \eqref{eq:condet2} implies that
$$
\dim(\mathcal{I}) \geq \ell-1.
$$ Thus \eqref{eq:condet1} holds.

Finally, note that \eqref{eq:condet2} is equivalent to the following statement: for any $\mathcal{I} \subseteq \mathbb{P}^{n-1}(\mathbb{F}_q)$ with $\dim(\langle \mathcal{I} \rangle)  \leq \ell-2$, we have $\sum_{{\bf w} \in \mathcal{I}}t_{{\bf w},S} \leq k-1$.
\end{proof}

\begin{lemma}\label{lem:thetincrink}
Let $2 \leq \ell \leq n$. Then $(\Theta_{n}^{(k,\ell)})_{k \geq \ell}$ is strictly increasing in $k$.
\end{lemma}

\begin{proof}
The proof follows the same line as that of \cite[Proposition 2.5 (1)]{AB24}. Let $S \subset \mathbb{B}_n$ be a $(k,\ell)$-orthogonal set of maximum size. Consider ${\bf v}\in \mathbb{B}_n \setminus S$. We claim that $S \cup \{{\bf v}\}$ is $(k+1,\ell)$-orthogonal. Indeed, let $X \subset S \cup \{{\bf v}\}$ be a subset of size $k+1$. Since $|S \cap X| \geq k$, the $(k,\ell)$-orthogonality of $S$ implies that $S \cap X$ contains an orthogonal set of size $\ell$. Thus, $S \cup \{{\bf v}\}$ is $(k+1,\ell)$-orthogonal.

Consequently, we have $\Theta_n^{(k,\ell)}=|S|<\Theta_n^{(k+1,\ell)}$.
\end{proof}
\subsection{Some Preliminaries on Geometry of Numbers in $\mathbb{F}_q((x^{-1}))^n$}
Suppose that $\mathcal{K}=\mathbb{F}_q((x^{-1}))$. An equivalent condition for the orthogonality of vectors in $\mathbb{B}_n$ arises from Hadamard's inequality (see, for example, \cite[Lemma 2.4]{RW}).

\begin{lemma}[Hadamard's Inequality]\label{lem:hadineq}
    Let $\mathbf{v}_1,\dots,\mathbf{v}_{\ell}\in \mathcal{K}^n$. Then, 
    \begin{equation}\label{eqn:Hadamard}\Vert \mathbf{v}_1\wedge \dots\wedge \mathbf{v}_{\ell}\Vert\leq \prod_{i=1}^{\ell}\Vert \mathbf{v}_i\Vert.\end{equation}
    Furthermore, equality holds in (\ref{eqn:Hadamard}) if and only if $\mathbf{v}_1,\dots,\mathbf{v}_{\ell}$ are orthogonal. 
\end{lemma}
One consequence of Hadamard's inequality and an analogue of Minkowski's second theorem \cite[Equations (24) and (25)]{Mah} is that the successive minima of a lattice $\Lambda \leq \mathcal{K}^n$ form an orthogonal $\mathbb{F}_q[x]$-basis of the lattice (see \cite{Ara} for details).

\begin{theorem}{\cite[Equations (24) and (25)]{Mah}}
\label{thm:Mink2nd}
    Let $\Lambda\leq \mathcal{K}^n$ be a lattice, and let $\mathbf{v}_1,\dots,\mathbf{v}_n\in \Lambda$ be a set of independent vectors satisfying $\Vert \mathbf{v}_i\Vert=\lambda_i(\Lambda)$. Then, 
    \begin{equation}
        \det(\Lambda)=\Vert \mathbf{v}_1\wedge \dots \wedge \mathbf{v}_n\Vert=\prod_{i=1}^n\Vert \mathbf{v}_i\Vert=\prod_{i=1}^n\lambda_i(\Lambda).
    \end{equation}
\end{theorem}
Moreover, as mentioned in the introduction, \cite{Ara} proves that all orthogonal $\mathbb{F}_q[x]$-bases of a lattice in $\mathcal{K}^n$ are composed of successive minima. 

\begin{theorem}{\cite[Lemma 1.18 and Theorem 1.21]{Ara}}
\label{thm:AraLattBase}
    Let $\Lambda \leq \mathcal{K}^n$ be a lattice. Then, an $\mathbb{F}_q[x]$-basis of $\Lambda$ is orthogonal if and only if it is a basis of successive minima. 
\end{theorem}
\section{Maximum orthogonal sets}\label{sec:maxorthset}
In this section, we investigate the size and structure of maximum orthogonal sets by proving Theorem \ref{thm:orthogset} and Theorem \ref{thm:structumaxsize}. We will use the same notation introduced in \S\ref{sec:prelim}. Additionally, let $\pi$ denote a fixed uniformizer of $\mathcal{K}$. For $1 \leq k \leq n$, let ${n \brack k}_q$ denote the number of $k$-dimensional subspaces of $\mathbb{F}_q^n$; more precisely
$$
{n \brack k}_q=\frac{\prod_{i=0}^{k-1}(q^n-q^i)}{\prod_{i=0}^{k-1}(q^k-q^i)}.
$$

\subsection{Proof of Theorem \ref{thm:orthogset}} 

     \subsubsection{Proof of Theorem \ref{thm:orthogset} \eqref{thm:orthogsetsizemaxTh}}\label{ssecproofofmainthm}
     We first prove the lower bound of \eqref{eq:boundalphakln}.\newline
        {\bf Lower bound of \eqref{eq:boundalphakln}.}
        Let $S_0 \subseteq \mathbb{B}_n$ be a set satisfying $\vert (\rho_n \circ \gamma_n)^{-1}({\bf z})\cap S_0\vert  = \left\lfloor \frac{(k-1)(q-1)}{q^{\ell-1}-1}\right\rfloor$ for all ${\bf z} \in \mathbb{P}^{n-1}(\mathbb{F}_q)$. Note that, for any $S' \subseteq S_0$ of size $k$, we have 
        $$\left|(\rho_n \circ \gamma_n)(S')\right| \geq \frac{q^{\ell-1}-1}{q-1}+1.$$ 
        Thus $|\gamma_n(S')|\geq q^{\ell-1}+1$ and hence $\dim(\langle \gamma_n(S') \rangle) \geq \ell$, which implies, by Remark \ref{rmk:subspee}, that $\dim(\langle (\rho_n \circ \gamma_n)(S') \rangle)\geq \ell-1$. Therefore, by Lemma \ref{lem:charoforthusingproj}, $S_0$ is $(k,\ell)$-orthogonal, and thus,
        $$
       \Theta_{n}^{(k,\ell)} \geq \vert S_0\vert  = \left\lfloor \frac{(k-1)(q-1)}{q^{\ell-1}-1}\right\rfloor \times\frac{q^n-1}{q-1}.
        $$ 
        \vskip 1mm
        \noindent
        {\bf-Upper bound of \eqref{eq:boundalphakln}.} Let $S \subseteq \mathbb{B}_n$ be a $(k,\ell)$-weakly orthogonal set of maximum size. By Lemma \ref{lem:rt}, we have
        \begin{equation}\label{eq:uppboundTheta}
        \begin{split}
                \vert S\vert = \sum_{{\bf w}\in \mathbb{P}^{n-1}(\mathbb{F}_q)}t_{{\bf w},S}=\frac{\sum_{\mathcal{I} \leq \mathbb{P}^{n-1}(\mathbb{F}_q), \dim( \mathcal{I} )=\ell-2}\sum_{{\bf w} \in  \mathcal{I}}t_{{\bf w},S}}{{n-1 \brack \ell-2}_q}\\
                \leq \frac{{n \brack \ell-1}_q(k-1)}{{n-1 \brack \ell-2}_q}= \frac{k-1}{q^{\ell-1}-1} \times (q^n-1),
        \end{split}
        \end{equation}
        so
        $$
        \Theta_{n}^{(k,\ell)}=\vert S\vert \leq \left\lfloor\frac{k-1}{q^{\ell-1}-1} \times (q^n-1)\right\rfloor.
        $$

\subsubsection{Proof of Theorem \ref{thm:orthogset} \eqref{thm:orthogsetthetadeltcomp}}
By definition, a subset of $\mathcal{K}^n$ of size $2$ is weakly orthogonal if and only if it is orthogonal. Thus, for $k\geq 2$, $(k,2)$-weakly orthogonal sets coincide with $(k,2)$-orthogonal sets. It follows immediately that $\Delta_n^{(k,2)}=\Theta_n^{(k,2)}$, so we may assume $\ell \geq 3$. 

First, let us address the case $(q,k,\ell)=(2,4,3)$. Since 
$$\left\lfloor\frac{k-1}{3}(2^n-1)\right\rfloor=2^n-1,$$ by Theorem \ref{thm:orthogset} \eqref{thm:orthogsetsizemaxTh}, we have $\Theta_{n}^{(4,3)}\leq 2^n-1$. Consider the set $$S'=\left\{\sum_{i\in I}\mathbf{e}_i \mid \emptyset\neq I\subseteq [n]\right\},$$ and let
$S=S'\cup\{\mathbf{e}_1+\pi\mathbf{e}_2\}$ which has size $2^n$, where $\pi$ denotes the fixed uniformizer of $\mathcal{K}$. Since $(\rho_n \circ \gamma_n)|_{S'}$ is injective, by Lemma \ref{lem:charoforthusingproj}, $S'$ is weakly orthogonal. Hence, for any $Y \subset S$ of size $4$, the subset $Y \cap S' \subset Y$ is weakly orthogonal and has size at least $3$. Therefore, $S$ is $(4,3)$-weakly orthogonal. Consequently
$$
\Delta_n^{(4,3)}\geq |S|=2^n > 2^n-1 \geq \Theta^{(4,3)}_{n}.
$$

Next, assume that $(q,k,\ell) \neq (2,4,3)$. By Theorem \ref{thm:orthogset} \eqref{thm:orthogsetsizemaxTh} and \cite[Theorem 1.10]{AB24}, to prove that $\Theta_n^{(k,\ell)}< \Delta_{n}^{(k,\ell)}$, it suffices to show that 
\begin{equation}\label{eq:tobeproved}
\left\lfloor\frac{k-1}{q^{\ell-1}-1}\times (q^n-1)\right\rfloor < \left\lfloor \frac{k-1}{\ell-1} \right\rfloor\frac{q^n-1}{q-1}.
\end{equation}
Suppose first that $(q,\ell) \in \{(3,3),(2,3),(2,4)\}$. We will check that \eqref{eq:tobeproved} is satisfied in each case.
\vskip 1mm
\noindent $\bullet$ {\bf Case 1: $(q,\ell)=(3,3)$.} Since $k \geq 3$, we have
$$
\frac{k-1}{8}< \frac{1}{2}\left\lfloor \frac{k-1}{2} \right\rfloor,
$$
so
$$
\left\lfloor\frac{k-1}{8}\times (3^n-1)\right\rfloor < \left\lfloor \frac{k-1}{2} \right\rfloor\frac{3^n-1}{2}.
$$
\vskip 1mm
\noindent $\bullet$ {\bf Case 2: $(q,\ell)=(2,3)$.} Since $k=3$ or $k \geq 5$, we have
$$
\frac{k-1}{3}< \left\lfloor \frac{k-1}{2} \right\rfloor,
$$
so
$$
\left\lfloor\frac{k-1}{3}\times (2^n-1)\right\rfloor < \left\lfloor \frac{k-1}{2} \right\rfloor\times (2^n-1).
$$
\vskip 1mm
\noindent $\bullet$ {\bf Case 3: $(q,\ell)=(2,4)$.} Since $k \geq 4$, we have
$$
\frac{k-1}{7}< \left\lfloor \frac{k-1}{3} \right\rfloor,
$$
so
$$
\left\lfloor\frac{k-1}{7}\times (2^n-1)\right\rfloor < \left\lfloor \frac{k-1}{3} \right\rfloor \times (2^n-1).
$$
\vskip 1mm
\noindent
Suppose now that $(q,\ell)\notin \{(3,3),(2,3),(2,4)\}$. Since $k \geq \ell$, we have
\begin{equation}\label{eq:lala1}
\frac{k-1}{\ell-1} \leq 2\left\lfloor \frac{k-1}{\ell-1} \right\rfloor.
\end{equation}
Moreover, it is easy to verify that
\begin{equation}\label{eq:lala2}
 2(\ell-1)<q^{\ell-2}+\dots+q+1=\frac{q^{\ell-1}-1}{q-1}.
\end{equation}
From \eqref{eq:lala1} and \eqref{eq:lala2}, we deduce
$$
\left\lfloor\frac{k-1}{q^{\ell-1}-1}\times (q^n-1)\right\rfloor \leq \frac{k-1}{q^{\ell-1}-1}\times (q^n-1)=\frac{k-1}{\ell-1}\times \frac{\ell-1}{q^{\ell-1}-1}\times(q^n-1)< \left\lfloor \frac{k-1}{\ell-1} \right\rfloor\frac{q^n-1}{q-1},
$$
and so \eqref{eq:tobeproved} holds. This completes the proof.

\subsubsection{Proof of Theorem \ref{thm:orthogset} \eqref{thm:orthogsettheind}}
By contrapositive, assume that $k > q^{\ell-1}$. Combining \cite[Theorem 1.23 (a)]{AB24} with Theorem \ref{thm:orthogset} \eqref{thm:orthogsetsizemaxTh}, we obtain 
${\rm Ind}_q(n,k,\ell)=q^n-1=\Theta_n^{(q^{\ell-1},\ell)}$. Now, by Lemma \ref{lem:thetincrink}, we conclude that
$$
{\rm Ind}_q(n,k,\ell)=\Theta_n^{(q^{\ell-1},\ell)}<\Theta_n^{(k,\ell)},
$$
as was to be proved.

\subsection{Proof of Theorem \ref{thm:structumaxsize}}
In this part, we prove Theorem \ref{thm:structumaxsize}. To do this, we need the following two lemmas.

\begin{lemma}\label{lem:weakorthsumeqval}
Let $f: \mathbb{P}^{n-1}(\mathbb{F}_q) \rightarrow \mathbb{Z}$ be a function, let $2 \leq \ell \leq \frac{q^n-1}{q-1}$ and let $r \in \mathbb{Z}$ be such that
$$
\sum_{{\bf z} \in \mathcal{I}} f({\bf z})=r,
$$ for all $\mathcal{I} \subset \mathbb{P}^{n-1}(\mathbb{F}_q)$ of size $\ell-1$. Then $f$ is constant, that is, $f=\frac{r}{\ell-1}$.
\end{lemma}

\begin{proof}
The proof is straightforward.
\end{proof}

\begin{lemma}\label{lem:orthsumtoeqval}
Let $f: \mathbb{P}^{n-1}(\mathbb{F}_q) \rightarrow \mathbb{Z}$ be a function, let $2 \leq \ell \leq n$ and let $r \in \mathbb{Z}$ be such that
$$
\sum_{{\bf z} \in \mathcal{I}} f({\bf z})=r,
$$ for all $\ell-2$-dimensional subspaces $\mathcal{I} \leq \mathbb{P}^{n-1}(\mathbb{F}_q)$. Then $f$ is constant, that is, $f=\frac{r(q-1)}{q^{\ell-1}-1}$.
\end{lemma}

\begin{proof}
The claim is trivial if $\ell=2$. Now, assume $\ell \geq 3$. Let ${\bf z}\neq {\bf z}' \in \mathbb{P}^{n-1}(\mathbb{F}_q)$. From the assumption, we have
\begin{equation}\label{eq:sumsubconzzpri}
\sum_{\mathcal{I}\,:\,  {\bf z} \in \mathcal{I},\,{\rm dim}(\mathcal{I})=\ell-2} \sum_{{\bf v} \in \mathcal{I}}f({\bf v})=r {n-1 \brack \ell-2}_q=\sum_{\mathcal{I}\,:\,{\bf z}' \in \mathcal{I},\,{\rm dim}(\mathcal{I})=\ell-2} \sum_{{\bf v} \in \mathcal{I}}f({\bf v}).
\end{equation}
Note that, for any ${\bf v} \in \mathbb{P}^{n-1}(\mathbb{F}_q)\setminus\{{\bf z},{\bf z}'\}$, the number $A_{{\bf v}}$ of $\ell-2$-dimensional subspaces $\mathcal{I} \leq \mathbb{P}^{n-1}(\mathbb{F}_q)$, such that $\mathbf{v},\mathbf{z}\in \mathcal{I}$ is equal to the number $B_{{\bf v}}$ of $\ell-2$-dimensional subspaces $\mathcal{I}'\leq \mathbb{P}^{n-1}(\mathbb{F}_q)$, such that $\mathbf{v},\mathbf{z}'\in \mathcal{I}'$. 
Moreover
\begin{align}
\sum_{\mathcal{I}\,:\,  {\bf z} \in \mathcal{I},\,{\rm dim}(\mathcal{I})=\ell-2} \sum_{{\bf v} \in \mathcal{I}}f({\bf v})&=\sum_{{\bf v} \in \mathbb{P}^{n-1}(\mathbb{F}_q)}f({\bf v})\left( \sum_{\mathcal{I}\,:\, {\bf v},{\bf z} \in \mathcal{I},\dim(\mathcal{I})=\ell-2} 1 \right)\nonumber\\
&=\sum_{\mathcal{I}\,:\,{\bf z} \in \mathcal{I},\,{\rm dim}(\mathcal{I})=\ell-2} f({\bf z})+  \left(\sum_{\mathcal{I}\,:\,{\bf z}, {\bf z}' \in \mathcal{I},\,{\rm dim}(\mathcal{I})=\ell-2} f({\bf z}')\right)\nonumber\\
&+ \sum_{{\bf v} \in \mathbb{P}^{n-1}(\mathbb{F}_q)\setminus\{{\bf z},{\bf z}'\}}f({\bf v}) \underbrace{\left(\sum_{\mathcal{I}\,:\,{\bf z}, {\bf v} \in \mathcal{I},\,{\rm dim}(\mathcal{I})=\ell-2} 1\right)}_{=A_{{\bf v}}},\nonumber
\end{align}
and
\begin{align}
\sum_{\mathcal{I}\,:\,  {\bf z}' \in \mathcal{I},\,{\rm dim}(\mathcal{I})=\ell-2} \sum_{{\bf v} \in \mathcal{I}}f({\bf v})&=\sum_{{\bf v} \in \mathbb{P}^{n-1}(\mathbb{F}_q)}f({\bf v})\left( \sum_{\mathcal{I}\,:\, {\bf v},{\bf z}' \in \mathcal{I},\dim(\mathcal{I})=\ell-2} 1 \right) \nonumber\\
&=\sum_{\mathcal{I}\,:\,{\bf z}' \in \mathcal{I},\,{\rm dim}(\mathcal{I})=\ell-2} f({\bf z}')+  \left(\sum_{\mathcal{I}\,:\,{\bf z}, {\bf z}' \in \mathcal{I},\,{\rm dim}(\mathcal{I})=\ell-2} f({\bf z})\right) \nonumber \\
&+ \sum_{{\bf v} \in \mathbb{P}^{n-1}(\mathbb{F}_q)\setminus\{{\bf z},{\bf z}'\}}f({\bf v}) \underbrace{\left(\sum_{\mathcal{I}\,:\,{\bf z}', {\bf v} \in \mathcal{I},\,{\rm dim}(\mathcal{I})=\ell-2} 1\right)}_{=B_{{\bf v}}}.\nonumber
\end{align}
Thus, from \eqref{eq:sumsubconzzpri}, we get

$$
\underbrace{\sum_{\mathcal{I}\,:\,{\bf z} \in \mathcal{I},\,{\rm dim}(\mathcal{I})=\ell-2} f({\bf z})}_{={n-1 \brack \ell-2 }_q f({\bf z})}+\underbrace{\sum_{\mathcal{I}\,:\,{\bf z}, {\bf z}' \in \mathcal{I},\,{\rm dim}(\mathcal{I})=\ell-2}f({\bf z}')}_{={n-2 \brack{\ell-3}}_q f({\bf z}')}=\underbrace{\sum_{\mathcal{I}\,:\,{\bf z}' \in \mathcal{I},\,{\rm dim}(\mathcal{I})=\ell-2} f({\bf z}')}_{={n-1 \brack{\ell-2}}_q f({\bf z}')}+\underbrace{\sum_{\mathcal{I}\,:\,{\bf z}, {\bf z}' \in \mathcal{I},\,{\rm dim}(\mathcal{I})=\ell-2}f({\bf z})}_{_{={n-2 \brack{\ell-3}}_q f({\bf z})}}.
$$ 
Hence, we obtain
$$
\left({n-1 \brack \ell-2 }_q-{n-2 \brack{\ell-3}}_q\right)(f({\bf z})-f({\bf z}'))=0.
$$ Since $\ell \leq n$, it follows that ${n-1 \brack \ell-2 }_q \neq {n-2 \brack{\ell-3}}_q$. Therefore, we conclude that $f({\bf z})=f({\bf z}')$. Consequently $f$ is constant on $\mathbb{P}^{n-1}(\mathbb{F}_q)$, that is, $f=\frac{r(q-1)}{q^{\ell-1}-1}$. This completes the proof. 
\end{proof}

\subsubsection{Proof of Theorem \eqref{thm:structumaxsize} \eqref{thm:structumaxsizeWOcase}: weakly orthogonal case} Let $2 \leq \ell \leq k$. Suppose that $k-1$ divides $\ell-1$.
Let $\mathcal{F}  \subset \mathbb{B}_n$ be a $(k,\ell)$-weakly orthogonal set of maximum size, that is, 
$$
|\mathcal{F}|=\frac{k-1}{\ell-1} \times \frac{q^n-1}{q-1}.
$$ By the proof of \cite[Theorem 1.10]{AB24}, we have $\sum_{{\bf w} \in \mathcal{I}}t_{{\bf w},\mathcal{F}}=k-1$ for all $\mathcal{I} \subset \mathbb{P}^{n-1}(\mathbb{F}_q)$ of size $\ell-1$. Now, applying Lemma \ref{lem:weakorthsumeqval} to the function $f=t_{\cdot,\mathcal{F}}$ and $r=k-1$, we obtain that $t_{\cdot, \mathcal{F}}$ is constant on $\mathbb{P}^{n-1}(\mathbb{F}_q)$, that is, $f=\frac{k-1}{\ell-1}$. This completes the proof. 

\subsubsection{Proof of Theorem \ref{thm:structumaxsize} \eqref{thm:structumaxsizeOrcase}: orthogonal case} Let $2 \leq \ell \leq k$. Suppose that $\frac{q^{\ell-1}-1}{q-1}$ divides $k-1$. Let $\mathcal{F} \subset \mathbb{B}_n$ be a $(k,\ell)$-orthogonal set of maximum size, that is, 
$$
|\mathcal{F}|=\frac{k-1}{q^{\ell-1}-1}\times (q^n-1).
$$
By Equation \eqref{eq:uppboundTheta}, we have $\sum_{{\bf w} \in  \mathcal{I}}t_{{\bf w},S}=k-1$ for all $\mathcal{I} \leq \mathbb{P}^{n-1}(\mathbb{F}_q)$ of dimension $\ell-2$. Now, applying Lemma \ref{lem:orthsumtoeqval} to the function $f=t_{\cdot,\mathcal{F}}$ and $r=k-1$, we obtain that $t_{\cdot, \mathcal{F}}$ is constant on $\mathbb{P}^{n-1}(\mathbb{F}_q)$, that is, $f=\frac{(k-1)(q-1)}{q^{\ell-1}-1}$. This completes the proof.

\section{Hadamard Matrices}
\label{sec:Hadamard}
In this section, we prove Theorem \ref{thm:hadamard} on Hadamard matrices. Throughout this section, let $\mathcal{K}=\mathbb{F}_q((x^{-1}))$ and $\mathcal{R}=\mathbb{F}_q[x]$. 

\subsection{Proof of Theorem \ref{thm:hadamard}}
Let $n\geq 1$, $\ell \geq 0$ and $m \leq \ell n$. 
\begin{proof}[Proof of Theorem \ref{thm:hadamard}\eqref{thm:hadamardcounting} and \eqref{thm:weakHadamarkCounting}]
Let $A \in M_{n \times n}(\mathcal{R}_{\leq \ell})$ be a matrix with nonzero columns, which are denoted by $\mathbf{v}_1,\dots,\mathbf{v}_n$. For every $1 \leq i \leq n$, let $\ell_i=\log_q(||{\bf v}_i||)\leq \ell$, and define ${\bf w}_i=x^{-\ell_i}{\bf v}_i \in \mathbb{B}_n$. Note that $A$ is $(\ell,q^m)$-Hadamard (resp., $(\ell,q^m)$-weakly Hadamard) if and only if $\ell_1+\dots+\ell_n=m$ and the tuple $({\bf w}_1,\dots,{\bf w}_n)$ is orthogonal (resp., weakly orthogonal); by Lemma \ref{lem:charoforthusingproj}, the latter means that $(\gamma_n(\mathbf{w}_1),\dots,\gamma_n(\mathbf{w}_n)) \in (\mathbb{F}_q^n)^n$ is independent (resp., weakly independent\footnote{that is, they are pairwise independent.}). Clearly, for an independent (resp., weakly independent) tuple  $({\bf a}_1,\dots,{\bf a}_n) \in (\mathbb{F}_q^n)^n$, we have $\left(\gamma(x^{-\ell_1}\mathbf{v}_1),\dots,\gamma(x^{-\ell_n}\mathbf{v}_n)\right)=(\mathbf{a}_1,\dots,\mathbf{a}_n)$, where $\Vert \mathbf{v}_i\Vert=q^{\ell_i}$, if and only if there exist $(\mathbf{u}_1,\dots,\mathbf{u}_n)\in \mathcal{R}^n$ with $\Vert \mathbf{u}_i\Vert\leq q^{\ell_i-1}$, such that $\mathbf{v}_i=x^{\ell_i}\mathbf{a}_i+\mathbf{u}_i$. Hence, there are exactly $q^{\ell_1} \times \dots \times q^{\ell_n}=q^m$ possible choices of $({\bf v}_1,\dots,{\bf v}_n)$ such that $\left(\gamma_n(x^{-\ell_1}\mathbf{v}_1),\dots,\gamma_n(x^{-\ell_n}\mathbf{v}_n)\right)=({\bf a}_1,\dots,{\bf a}_n)$ and $\Vert \mathbf{v}_i\Vert=q^{\ell_i}$ for every $i=1,\dots,n$. 

Note that the number of linearly independent tuples $(\mathbf{a}_1,\dots,\mathbf{a}_n)\in (\mathbb{F}_q^n)^n$ is equal to $\prod_{i=0}^{n-1}(q^n-q^i)$, since this is the cardinality of $\operatorname{GL}_n(\mathbb{F}_q)$ \cite[page 35, equation (2)]{DummitFoote2004}. Moreover, by \cite[Lemma 3.7(1)]{AB24}, the set $\left\{(\rho_n\circ \gamma_n)(x^{-\ell_i}\mathbf{v}_i):i\in[n]\right\}=\left\{\rho_n(\mathbf{a}_i):i\in [n]\right\}$ has cardinality $n$. Note that $\vert \mathbb{P}^{n-1}(\mathbb{F}_q)\vert=\frac{q^n-1}{q-1}$, and therefore, the number tuples $(\mathbf{a}_1,\dots,\mathbf{a}_n)\in (\mathbb{P}^{n-1}(\mathbb{F}_q))^n$ with distinct coordinates is $\begin{pmatrix}
    \frac{q^{n-1}}{q-1}\\
    n
\end{pmatrix}n!$. Hence, the number of weakly independent tuples $(\mathbf{a}_1,\dots,\mathbf{a}_n)\in (\mathbb{F}_q^n)^n$ is $(q-1)^nn!\begin{pmatrix}
    \frac{q^{n}-1}{q-1}\\
    n
\end{pmatrix}$. Since $({\bf v}_1,\dots,{\bf v}_n)$ is uniquely determined by $(\ell_1,\dots,\ell_n)$ and $({\bf w}_1,\dots,{\bf w}_n)$, we obtain
$$
|\mathcal{H}_n(\ell,q^m)|=q^m \prod_{i=0}^{n-1}(q^n-q^i)\left\vert\left\lbrace(\ell_1,\dots,\ell_n) \in \{0,\dots,\ell\}^n \Bigm\vert \sum_{i=1}^n \ell_i=m\right\rbrace\right\vert,
$$
and
$$
|\mathcal{W}_n(\ell,q^m)|=q^m (q-1)^n n! \binom{\frac{q^n-1}{q-1}}{n} \left\vert\left\lbrace(\ell_1,\dots,\ell_n) \in \{0,\dots,\ell\}^n \Bigm\vert \sum_{i=1}^n \ell_i=m\right\rbrace\right\vert.
$$ 
\end{proof}

\subsection{Proof of Theorem \ref{thm:hadamardKroneckerprod}} 
Let $n_1,n_2 \geq 1$, $\ell_1,\ell_2 \geq 0$ and $m_1 \leq \ell_1 n_1,m_2 \leq \ell_2 n_2$.

\begin{proof}[Proof of Theorem \ref{thm:hadamardKroneckerprod} for Hadamard Matrices]
 Let $A \in M_{n_1\times n_1}(\mathcal{R}_{\leq \ell_1})$ (resp., $B \in M_{n_2\times n_2}(\mathcal{R}_{\leq \ell_2})$) be an $(\ell_1,q^{m_1})$-Hadamard (resp., $(\ell_2,q^{m_2})$-Hadamard) matrix. Clearly, we see that $A \otimes B \in M_{n_1n_2,n_1n_2}(\mathcal{R}_{\leq \ell_1+\ell_2})$. Now, denote the columns of $A$ (resp., $B$) by $A_1,\dots,A_{n_1}$ (resp., $B_1,\dots,B_{n_2}$). By definition, the entries of $A \otimes B$ are given by: 
 $$(A\otimes B)_{n_2(r-1)+i,n_2(s-1)+j}=a_{rs}b_{ij},$$ for all 
 $1 \leq r,s \leq n_1$ and $1 \leq i,j \leq n_2$.
Thus, the norm of the $n_2(s-1)+j$-th column of $A\otimes B$ is $\max_{1 \leq r \leq n_1,1 \leq i \leq n_2}\vert a_{rs}b_{ij}\vert=\Vert A_s\Vert\times \Vert B_j\Vert$. Therefore, we obtain
    \begin{equation}
    \begin{split}
        \prod_{j=1}^{n_2}\prod_{s=1}^{n_1}\Vert (A\otimes B)_{n_2(s-1)+j}\Vert&=\prod_{j=1}^{n_2}\prod_{s=1}^{n_1}(\Vert A_s\Vert\times \Vert B_j\Vert)\\
        &=\prod_{j=1}^{n_2}\Vert B_j\Vert^{n_1}\left(\prod_{s=1}^{n_1}\Vert A_s\Vert\right)\\
        &=\left(\prod_{s=1}^{n_1}\Vert A_s\Vert^{n_2}\right)\times \left(\prod_{j=1}^{n_2}\Vert B_j\Vert^{n_1}\right)\\
        &=\vert \det(A)\vert^{n_2}\vert \det(B)\vert^{n_1}\quad\textrm{by the Hadamard property}\\
        &=\vert \det(A\otimes B)\vert.
    \end{split}
    \end{equation}
    By Lemma \ref{lem:hadineq}, the columns vectors of $A \otimes B$ are orthogonal. Consequently, $A\otimes B$ is $(\ell_1+\ell_2,q^{m_1n_2+m_2n_1})$-Hadamard.
\end{proof}

\begin{proof}[Proof of Theorem \ref{thm:hadamardKroneckerprod} for Weakly Hadamard Matrices]
 Let $A \in M_{n_1\times n_1}(\mathcal{R}_{\leq \ell_1})$ (resp., $B \in M_{n_2\times n_2}(\mathcal{R}_{\leq \ell_2})$) be an $(\ell_1,q^{m_1})$-weakly Hadamard (resp., $(\ell_2,q^{m_2})$-weakly Hadamard) matrix.  Clearly, we see that $A \otimes B \in M_{n_1n_2,n_1n_2}(\mathcal{R}_{\leq \ell_1+\ell_2})$. Moreover, by \cite[Lemma 3.2]{AB24}, 
 \begin{equation}\label{eq:indepengammee}
 \gamma_{n_1}\left(\frac{A_1}{x^{\log_q\Vert A_1\Vert}}\right),\dots,\gamma_{n_1}\left(\frac{A_{n_1}}{x^{\log_q\Vert A_{n_1}\Vert}}\right)\,\,\,  \textrm{are pairwise independent over}\,\, \mathbb{F}_q.\end{equation}
 Similarly, 
 \begin{equation}\label{eq:indepengammee2}
 \gamma_{n_2}\left(\frac{B_1}{x^{\log_q\Vert B_1\Vert}}\right),\dots,\gamma_{n_2}\left(\frac{B_{n_2}}{x^{\log_q\Vert B_{n_2}\Vert}}\right)\,\,\,  \textrm{are pairwise independent over}\,\, \mathbb{F}_q.
 \end{equation}
 Again, by \cite[Lemma 3.2]{AB24}, to prove that the columns of $A\otimes B$ are weakly orthogonal, we need to prove that the vectors 
\begin{equation}\label{eq:indepkroneck}
\left\{\gamma_{n_1n_2}\left(\frac{(A\otimes B)_{n_2(s-1)+j}}{x^{\log_q\Vert (A\otimes B)_{n_2(s-1)+j}\Vert}}\right)\right\}_{1\leq s\leq n_1,1\leq j\leq n_2}\,\,\,\textrm{are pairwise independent over}\,\,\mathbb{F}_q.
\end{equation} Since $\Vert (A\otimes B)_{n_2(s-1)+j}\Vert=\Vert A_s\Vert \times \Vert B_j\Vert$, we have 
 $$\log_q\Vert (A\otimes B)_{n_2(s-1)+j}\Vert=\log_q\Vert A_s\Vert+\log_q\Vert B_j\Vert,$$ for all $1 \leq s \leq n_1$ and $1 \leq j \leq n_2$. Thus, 
$$\left(\frac{(A\otimes B)_{n_2(s-1)+j}}{x^{\log_q\Vert A\otimes B\Vert_{n_2(s-1)+j}\Vert}}\right)_{n_2(r-1)+i}=\frac{a_{rs}}{x^{\log_q\Vert A_s\Vert}}\frac{b_{ij}}{x^{\log_q\Vert B_j\Vert}},$$
for all 
 $1 \leq r,s \leq n_1$ and $1 \leq i,j \leq n_2$.
Hence 
\begin{equation}\label{eq:gamman1n2coeff}
\left(\gamma_{n_1n_2}\left(\frac{(A\otimes B)_{n_2(s-1)+j}}{x^{\log_q\Vert (A\otimes B)_{n_2(s-1)+j}\Vert}}\right)\right)_{n_2(r-1)+i}=\gamma_1\left(\frac{a_{rs}}{x^{\log_q\Vert A_s\Vert}}\right)\gamma_1\left(\frac{b_{ij}}{x^{\log_q\Vert B_j\Vert}}\right),
\end{equation}
for all 
 $1 \leq r,s \leq n_1$ and $1 \leq i,j \leq n_2$.

Let $1\leq s,s'\leq n_2$ and $1\leq j,j'\leq n_1$ be such that  
$$c\gamma_{n_1n_2}\left(\frac{(A\otimes B)_{n_2(s-1)+j}}{x^{\log_q\Vert (A\otimes B)_{n_2(s-1)+j}\Vert}}\right)=\gamma_{n_1n_2}\left(\frac{(A\otimes B)_{n_2(s'-1)+j'}}{x^{\log_q\Vert (A\otimes B)_{n_2(s'-1)+j'}\Vert}}\right),
$$ for some $c\in \mathbb{F}_q$. Then, by \eqref{eq:gamman1n2coeff}, we have 
\begin{equation}
\gamma_1\left(\frac{a_{rs'}}{x^{\log_q\Vert A_{s'}\Vert}}\right)\gamma_1\left(\frac{b_{i,j'}}{x^{\log_q\Vert B_{j'}\Vert}}\right)=c\gamma_1\left(\frac{a_{rs}}{x^{\log_q\Vert A_s\Vert}}\right)\gamma_1\left(\frac{b_{ij}}{x^{\log_q\Vert B_j\Vert}}\right),
\end{equation}
for all $1\leq r\leq n_1$ and $1\leq i\leq n_2$.
Since $\gamma_{n_2}\left(\frac{B_{j'}}{x^{\log_q\Vert B_{j'}\Vert}}\right) \neq {\bf 0}$, there exists $1 \leq i_0 \leq n_2$ such that $\gamma_{1}\left(\frac{b_{i_0j'}}{x^{\log_q\Vert B_{j'}\Vert}}\right) \neq 0$. Thus, 
\begin{equation}
    \gamma_1\left(\frac{a_{rs'}}{x^{\log_q\Vert A_{s'}\Vert}}\right)=\overbrace{c\frac{\gamma_1\left(\frac{b_{i_0j}}{x^{\log_q\Vert B_j\Vert}}\right)}{\gamma_1\left(\frac{b_{i_0j'}}{x^{\log_q\Vert B_{j'}\Vert}}\right)}}^{=C}\gamma_1\left(\frac{a_{rs}}{x^{\log_q\Vert A_s\Vert}}\right),
\end{equation}
for all $1 \leq r \leq n_1$.
Hence,
$$\gamma_{n_1}\left(\frac{A_{s'}}{x^{\log_q\Vert A_{s'}\Vert}}\right)=C\gamma_{n_1}\left(\frac{A_{s}}{x^{\log_q\Vert A_{s}\Vert}}\right).$$ By \eqref{eq:indepengammee}, we deduce $s=s'$. Using with the same argument with \eqref{eq:indepengammee2}, we obtain $j=j'$. Therefore, \eqref{eq:indepkroneck} holds, so the columns of $A\otimes B$ are weakly orthogonal. 

Noting that $\vert \det(A\otimes B)\vert=q^{n_1m_2+n_2m_1}$, we conclude that $A\otimes B$ is $(\ell_1+\ell_2,q^{n_1m_2+n_2m_1})$-weakly Hadamard. 
\end{proof}

\section{On the sublattices of $\mathbb{F}_q((x^{-1}))^n$}
\label{sec:SubLattice}

In this section, we prove Theorem \ref{thm:AraBase}, Theorem \ref{thm:numOrtbases}, Theorem \ref{thm:NumG_n,k,q(mu,s)} and Theorem \ref{thm:PrimLattCount},  which address counting problems for sublattices. For this purpose, let $\mathcal{K}=\mathbb{F}_{q}((x^{-1}))$ and $\mathcal{R}=\mathbb{F}_q[x]$. 
\subsection{Proof of Proposition \ref{thm:AraBase}} 
Recall that Proposition \ref{thm:AraBase} is an analogue of Theorem \ref{thm:AraLattBase} for $k$-sublattices. To proceed with its proof, we first prove that analogous versions of Theorem \ref{thm:Mink2nd} and Theorem \ref{thm:AraBase} hold for $k$-sublattices.
\begin{lemma}
\label{lem:k-subOrt}
    Let $\Lambda\leq \mathcal{K}^n$ be a $k$-sublattice, and let $\mathbf{v}_1,\dots,\mathbf{v}_k\in \Lambda$ be a set of linearly independent vectors such that $\lambda_i(\Lambda)=\Vert \mathbf{v}_i\Vert$. Then, $\mathbf{v}_1,\dots,\mathbf{v}_k$ are orthogonal.
\end{lemma}

\begin{proof}
    Firstly, note that if $\mathbf{u}_1,\dots,\mathbf{u}_{\ell}$ are orthogonal, then, for any non-zero scalars $c_1,\dots,c_{\ell}\in \mathcal{K}\setminus\{0\}$, the vectors $c_1\mathbf{u}_1,\dots,c_{\ell}\mathbf{u}_{\ell}$ are orthogonal. Hence, by replacing $\mathbf{v}_i$ with $\frac{\mathbf{v}_i}{x^{\log_q\lambda_i(\Lambda)}}$, it suffices to assume that $\mathbf{v}_1,\dots,\mathbf{v}_k\in \mathbb{B}_n$. By Proposition \ref{prop:OrtIndConn}, if $\mathbf{v}_1,\dots,\mathbf{v}_k$ are not orthogonal, then, $\gamma_n(\mathbf{v}_1),\dots,\gamma_n(\mathbf{v}_k)\in \mathbb{F}_q^n$ are not linearly independent. Hence, there exist $a_1,\dots,a_k\in \mathbb{F}_q$ which are not all zero, such that $\sum_{i=1}^ka_i\gamma_n(\mathbf{v}_i)=0$. As a consequence, $\gamma\left(\sum_{i=1}^ka_i\mathbf{v}_i\right)=0$, so that $\left\Vert\sum_{i=1}^k a_i\mathbf{v}_i\right\Vert<1=\lambda_1(\Lambda)$. On the other hand, since $a_i\in \mathbb{F}_q\subseteq \mathcal{R}$, the non-zero vector $\mathbf{v}=\sum_{i=1}^ka_i\mathbf{v}_i$ satisfies $\mathbf{v}\in\Lambda$. Therefore, $\left\Vert\mathbf{v}\right\Vert\geq \lambda_1(\Lambda)=1$, which is a contradiction. Thus, $\mathbf{v}_1,\dots,\mathbf{v}_k$ are orthogonal. 
\end{proof}

\begin{lemma}
\label{rmk:SuccMin}
    Let $\Lambda \leq \mathcal{K}^n$ be a $k$-sublattice, and let $\{\mathbf{v}_1,\dots,\mathbf{v}_k\} \subset \Lambda$ be a basis of successive minima for $\Lambda$. Suppose that $\mathbf{v}\in \Lambda$ satisfies $\Vert \mathbf{v}\Vert<\lambda_i(\Lambda)$, for some $2 \leq i \leq k$. Then, $\mathbf{v}$ is an $\mathcal{R}$-linear combination of $\mathbf{v}_1,\dots,\mathbf{v}_{i-1}$.
\end{lemma}

\begin{proof}
    We prove this claim inductively. Let $\mathbf{v}\in \Lambda$ be a vector satisfying $\Vert \mathbf{v}\Vert<\lambda_i(\Lambda)$. If $\mathbf{v}\notin \operatorname{span}_{\mathcal{K}}\left\{\mathbf{v}_1,\dots,\mathbf{v}_{i-1}\right\}$, then, by the definition of successive minima, $\Vert \mathbf{v}\Vert\geq \lambda_i(\Lambda)$. As a consequence, $\mathbf{v}\in \operatorname{span}_{\mathcal{K}}\left\{\mathbf{v}_1,\dots,\mathbf{v}_{i-1}\right\}$. Hence, there exist $a_1,\dots,a_{i-1}\in \mathcal{R}$ and $\alpha_1,\dots,\alpha_{i-1}\in \mathfrak{m}$, such that $\mathbf{v}=\sum_{j=1}^{i-1}(a_j+\alpha_j)\mathbf{v}_j$. Note that $\sum_{j=1}^{i-1}a_j\mathbf{v}_j\in \operatorname{span}_{\mathcal{R}}\left\{\mathbf{v}_1,\dots,\mathbf{v}_{i-1}\right\}$, and therefore, it suffices to prove that $\sum_{j=1}^{i-1}\alpha_j\mathbf{v}_j\in \operatorname{span}_{\mathcal{R}}\left\{\mathbf{v}_1,\dots,\mathbf{v}_{i-1}\right\}$. Firstly, note that by Lemma \ref{lem:k-subOrt},
    \begin{equation}
        \left\Vert\sum_{j=1}^{i-1}\alpha_j\mathbf{v}_j\right\Vert=\max_{j=1,\dots,i-1}\vert \alpha_j\vert\lambda_j(\Lambda)\leq \frac{1}{q}\lambda_{i-1}(\Lambda).
    \end{equation}
    As a consequence, $\alpha_{i-1}=0$. We prove inductively that $\alpha_j=0$ for every $j=1,\dots,i-1$. Assume that we have proved that $\alpha_{j+1}=\dots=\alpha_{i-1}=0$, and we now prove that $\alpha_j=0$. Note that in this case, $\sum_{\ell=1}^{i-1}\alpha_{\ell}\mathbf{v}_{\ell}=\sum_{\ell=1}^j\alpha_{\ell}\mathbf{v}_{\ell}$. Moreover, by Lemma \ref{lem:k-subOrt},
    \begin{equation}
        \left\Vert\sum_{\ell=1}^j\alpha_{\ell}\mathbf{v}_{\ell}\right\Vert\leq \max_{\ell=1,\dots,j}\vert \alpha_{\ell}\vert\lambda_{\ell}(\Lambda)\leq \frac{1}{q}\lambda_{j}(\Lambda).
    \end{equation}
    Hence, $\alpha_j=0$, so that by continuing inductively, we obtain that $\sum_{j=1}^{i-1}\alpha_j\mathbf{v}_j=0$. Thus, $\mathbf{v}=\sum_{j=1}^{i-1}a_j\mathbf{v}_j\in \operatorname{span}_{\mathcal{R}}\left\{\mathbf{v}_1,\dots,\mathbf{v}_{i-1}\right\}$.
\end{proof}

\begin{lemma}\label{rm:succindimplybasss}
Suppose ${\bf v}_1,\dots,{\bf v}_k \in \Lambda$ are $\mathcal{K}$-linearly  independent and satisfy $\Vert {\bf v}_i \Vert=\lambda_i(\Lambda)$ for all $1 \leq i \leq k$. Then $\{{\bf v}_1,\dots,{\bf v}_k\}$ forms an $\mathcal{R}$-basis of successive minima for $\Lambda$. 
\end{lemma}

\begin{proof}
    This proof follows by an argument similar to that used in the proof of \cite[Lemma 1.18]{Ara}. Since $\mathbf{v}_1,\dots,\mathbf{v}_k\in \Lambda$ are $\mathcal{K}$-linearly independent, they form a $k$-sublattice in $\Lambda$. Thus, it suffices to prove that every $\mathbf{v}\in \Lambda$ can be written as a linear combination over $\mathcal{R}$ of $\mathbf{v}_1,\dots,\mathbf{v}_k$. Since $\Lambda$ is a $k$-sublattice, then, $\Lambda \subseteq \operatorname{span}_{\mathcal{K}}\{\mathbf{v}_1,\dots,\mathbf{v}_k\}$. Thus, for every $\mathbf{v}\in \Lambda$, there exist $a_1,\dots,a_k\in \mathcal{R}$ and $\alpha_1,\dots,\alpha_k\in \mathfrak{m}$, such that $\mathbf{v}=\sum_{i=1}^k(a_i+\alpha_i)\mathbf{v}_i$. 
    
    We now prove inductively that $\alpha_i=0$ for every $i$. For $i=k$, note that by Lemma \ref{lem:k-subOrt},
    \begin{equation}
        \left\Vert\sum_{i=1}^k\alpha_i\mathbf{v}_i\right\Vert=\max_{i=1,\dots,k}\vert \alpha_i\vert\lambda_i(\Lambda)\leq \frac{1}{q}\lambda_k(\Lambda).
    \end{equation}
    Thus, by Lemma \ref{rmk:SuccMin}, $\alpha_k=0$. Assume that $\alpha_{j+1}=\dots=\alpha_k=0$. Then, $\sum_{i=1}^j\alpha_i\mathbf{v}_i\in \Lambda$ and
    \begin{equation}
        \left\Vert\sum_{i=1}^j\alpha_i\mathbf{v}_i\right\Vert=\max_{i=1,\dots,j}\vert \alpha_i\vert\lambda_i(\Lambda)\leq \frac{1}{q}\lambda_j(\Lambda),
    \end{equation}
    so that by Lemma \ref{rmk:SuccMin}, $\alpha_j=0$. Therefore, $\alpha_i=0$ for every $i=1,\dots,k$, so that $\mathbf{v}=\sum_{i=1}^ka_i\mathbf{v}_i\in \operatorname{span}_{\mathcal{R}}\{\mathbf{v}_1,\dots,\mathbf{v}_k\}$. Hence, $\Lambda=\operatorname{span}_{\mathcal{R}}\left\{\mathbf{v}_1,\dots,\mathbf{v}_k\right\}$.
\end{proof}

\begin{proof}[Proof of Proposition \ref{thm:AraBase}]
Lemma \ref{lem:k-subOrt} and Lemma \ref{rm:succindimplybasss} imply that if $\{\mathbf{u}_1,\dots,\mathbf{u}_k\}\subseteq \Lambda$ is a set of successive minima, then $\Lambda=\operatorname{span}_{\mathcal{R}}\left\{\mathbf{u}_1,\dots,\mathbf{u}_k\right\}$. Thus, it suffices to prove that an orthogonal $\mathcal{R}$ basis of a $k$-sublattice $\Lambda$ consists of successive minima. Let $\Lambda$ be a $k$-sublattice, and let $\{\mathbf{v}_1,\dots,\mathbf{v}_k\}$ be an orthogonal $\mathcal{R}$ basis of $\Lambda$ satisfying $\Vert \mathbf{v}_1\Vert\leq \dots\leq \Vert \mathbf{v}_k\Vert$. In a similar manner as done in \cite[Theorem 1.21]{Ara}, if $\{\mathbf{v}_1,\dots,\mathbf{v}_k\}$ is not a set of successive minima, then there exists $i$ such that $\Vert \mathbf{v}_{i}\Vert>\lambda_{i}(\Lambda)$. Let $i_0$ be the minimal $i$ such that $\Vert \mathbf{v}_i\Vert>\lambda_i(\Lambda)$. Since $\{\mathbf{v}_1,\dots,\mathbf{v}_k\}$ is an orthogonal $\mathcal{R}$ basis of $\Lambda$, then, by Lemma \ref{lem:k-subOrt} and Lemma \ref{rm:succindimplybasss}, 
\begin{equation}
    \label{eqn:det(Lambda)=Ort}
    \prod_{i=1}^k\lambda_i(\Lambda)=\vert \det(\Lambda)\vert=\Vert\mathbf{v}_1\wedge\dots\wedge \mathbf{v}_k\Vert=\prod_{i=1}^k\Vert \mathbf{v}_i\Vert.
\end{equation}
Since $\Vert \mathbf{v}_{i_0}\Vert>\lambda_{i_0}(\Lambda)$, equation \eqref{eqn:det(Lambda)=Ort} implies that there exists $i_1>i_0$, such that $\Vert \mathbf{v}_{i_1}\Vert<\lambda_{i_1}(\Lambda)$. Thus, by Lemma \ref{rmk:SuccMin}, for every $i\leq i_1$, we have $\mathbf{v}_{i}\in \operatorname{span}_{\mathcal{R}}\left\{\mathbf{u}_1,\dots,\mathbf{u}_{i_1-1}\right\}$, where $\{\mathbf{u}_1,\dots,\mathbf{u}_k\}$ are successive minima of $\Lambda$. Thus, 
$$\dim\operatorname{span}\{\mathbf{v}_1,\dots,\mathbf{v}_{i_1}\}\leq \dim\operatorname{span}\{\mathbf{u}_1,\dots,\mathbf{u}_{i_1-1}\}\leq i_1-1,$$
which is a contradiction. As a consequence, $\left\{\mathbf{v}_1,\dots,\mathbf{v}_k\right\}$ is a set of successive minima of $\Lambda$.
\end{proof}
\subsection{Proof of Theorem \ref{thm:numOrtbases}}
In this part, we provide a formula for the number of orthogonal bases of sublattices of $\mathcal{K}^n$. To do so, by Proposition \ref{thm:AraBase}, it suffices to find necessary and sufficient conditions for a given basis of a sublattice to be a basis of successive minima. Let $\Lambda$ be a $k$-dimensional $(\boldsymbol{\mu},\boldsymbol{s})$-sublattice (see Definition \ref{def:ussublatt}), where $\boldsymbol{\mu}=(\mu_1,\dots,\mu_m)\in (q^\mathbb{Z})^m$ and ${\bf s} \in \mathbb{N}^m$ are as in Definition \ref{def:ussublatt}. For brevity, let $\lambda_i=\lambda_i(\Lambda)$ for all $1 \leq i \leq k$. Fix an ordered basis $\{\mathbf{v}_1,\dots,\mathbf{v}_k\}$ of successive minima for $\Lambda$. For $1 \leq \ell \leq m$, let $$I_{\ell}=\{j\in \{1,\dots,k\} \mid \lambda_j=\mu_{\ell}\}=\left\{\sum_{i=1}^{\ell-1}s_i+1,\dots,\sum_{i=1}^{\ell-1}s_i+s_{\ell}\right\}.$$

For vectors
\begin{equation}\label{eqn:u_1...u_n}\mathbf{u}_1=\sum_{i=1}^ka_{i,1}\mathbf{v}_i \in \Lambda,\dots,\mathbf{u}_k=\sum_{i=1}^ka_{i,k}\mathbf{v}_i \in \Lambda,\,\, (a_{i,j}\in \mathcal{R}),
\end{equation} 
consider the following conditions:
\begin{enumerate}
\item \label{cond:I_ellCond} $\vert a_{i,j}\vert\leq 1$ (that is, $a_{i,j}\in \mathbb{F}_q$), for $i,j\in I_{\ell}$ with $1 \leq \ell \leq m$;
\item \label{cond:l'<lCond}$\vert a_{i,j}\vert\leq \frac{\mu_{\ell}}{\mu_{\ell'}}$, for all $i \in I_{\ell'},j \in I_\ell$ with $1 \leq \ell'<\ell \leq m$;
\item \label{cond:l'>lCond}$a_{i,j}=0$, for all $i\in I_{\ell'}, j \in I_\ell$ with $1 \leq \ell< \ell' \leq m$;
 \item \label{cond:invertAl} the matrix $A_\ell=(a_{i,j})_{i,j\in I_{\ell}}$ is invertible over $\mathbb{F}_q$, for all $1 \leq \ell \leq m$.
\end{enumerate}

\begin{remark} \label{lem:a_sigma(p),p!=0}
The condition \eqref{cond:invertAl} implies that for every $1 \leq \ell \leq m$, there exists a permutation $\sigma_{\ell}$ on $I_{\ell}$ such that $a_{\sigma_\ell(j),j}\neq 0$ for all $j\in I_{\ell}$.
\end{remark}
We first prove that that $\{{\bf u}_1,\dots,{\bf u}_k\}$ is a basis of successive minima if and only if \eqref{cond:I_ellCond}, \eqref{cond:l'<lCond}, \eqref{cond:l'>lCond} and \eqref{cond:invertAl} are satisfied.

\begin{proposition}\label{prop:condnecfororth}
Let $\Lambda\in G_{n,k,q}(\boldsymbol{\mu},\boldsymbol{s})$ and fix an ordered basis of successive minima for $\Lambda$, $(\mathbf{v}_1,\dots,\mathbf{v}_k)$. Suppose that $( {\bf u}_1,\dots,{\bf u}_k)$ is a basis of successive minima for $\Lambda$ of the form \eqref{eqn:u_1...u_n}. Then the set $(\mathbf{u}_1,\dots,\mathbf{u}_k)$ satisfies conditions \eqref{cond:I_ellCond}, \eqref{cond:l'<lCond}, \eqref{cond:l'>lCond} and \eqref{cond:invertAl}.
\end{proposition}

\begin{proof}
First we prove that \eqref{cond:I_ellCond}, \eqref{cond:l'<lCond} and \eqref{cond:l'>lCond} are satisfied. For that, note that, by the orthogonality of $\mathbf{v}_1,\dots,\mathbf{v}_k$, we have
\begin{equation}
 \lambda_j=\Vert {\bf u}_j\Vert=\left\Vert\sum_{i=1}^ka_{i,j}\mathbf{v}_i\right\Vert=\max_{1 \leq i \leq k}\vert a_{i,j}\vert\times \Vert \mathbf{v}_i\Vert= \max_{1 \leq i \leq k} \lambda_i \vert a_{i,j}\vert ,
\end{equation}
for all $1 \leq j \leq k$. Therefore, $\vert a_{i,j}\vert\leq \frac{\lambda_j}{\lambda_i}$ for all $1 \leq i,j \leq k$. It follows immediately that \eqref{cond:I_ellCond}, \eqref{cond:l'<lCond} and \eqref{cond:l'>lCond} are satisfied.

Next we prove \eqref{cond:invertAl}. Assume, on the contrary, that $A_{\ell}$ is not invertible for some $1 \leq \ell \leq m$. Then, there exist $\{t_j\}_{j\in I_\ell}\subseteq \mathbb{F}_q$, not all zero, such that
    \begin{equation}\label{eq:sumaijtiszero}
        \sum_{j\in I_{\ell}}a_{i,j}t_j=0,
    \end{equation} for all $i\in I_{\ell}$.
    Thus, we have
    \begin{equation}
    \label{eqn:t_ju_jSum}
    \begin{split}
        \sum_{j\in I_{\ell}}t_j\mathbf{u}_j&=\sum_{j\in I_{\ell}}t_j\left(\sum_{i\in I_{\ell}}a_{i,j}\mathbf{v}_i+\sum_{\ell'<\ell}\sum_{i\in I_{\ell'}}a_{i,j}\mathbf{v}_i\right) \quad \quad \quad \quad \textrm{by \eqref{cond:l'>lCond}}\\
        &=\sum_{i\in I_{\ell}}\left(\sum_{j\in I_{\ell}}a_{i,j}t_j \right) \mathbf{v}_i+\sum_{j\in I_{\ell}}\sum_{\ell'<\ell}\sum_{i\in I_{\ell}'}t_ja_{i,j}\mathbf{v}_i \\
        &=\sum_{j\in I_{\ell}}\sum_{\ell'<\ell}\sum_{i\in I_{\ell}'}t_ja_{i,j}\mathbf{v}_i\in \operatorname{span}_{\mathcal{R}}\left\{\mathbf{v}_i \mid i\in \bigcup_{\ell'<\ell}I_{\ell'}\right\}\quad \textrm{by \eqref{eq:sumaijtiszero}}.
    \end{split}
    \end{equation}
    To derive a contradiction, observe that, for each $1 \leq i\leq \sum_{\ell'\leq \ell-1}s_{\ell'}$, we have 
    $$\Vert \mathbf{v}_i\Vert\leq \mu_{\ell-1}=\lambda_{\sum_{i=1}^{\ell-1}s_i}<\lambda_{\sum_{i=1}^{\ell-1}s_i+1},$$ so that, by Lemma \ref{rmk:SuccMin}, $\mathbf{v}_i$ is a linear combination of $\mathbf{u}_1,\dots,\mathbf{u}_{\sum_{i=1}^{\ell-1}s_i}$ over $\mathcal{R}$. Combining this with (\ref{eqn:t_ju_jSum}), we get $0 \neq \sum_{j\in I_{\ell}}t_j\mathbf{u}_j\in \operatorname{span}\left\{\mathbf{u}_1,\dots,\mathbf{u}_{\sum_{i=1}^{\ell-1}s_i}\right\}$, which contradicts the linear independence of ${\bf u}_1,\dots,{\bf u}_k$. Consequently, $A_\ell$ is invertible for all $1 \leq \ell \leq m$, so \eqref{cond:invertAl} is satisfied.
\end{proof}

\begin{proposition}\label{prop:condsuffororth}
Suppose that $({\bf u}_1,\dots,{\bf u}_k)$ satisfy \eqref{cond:I_ellCond}, \eqref{cond:l'<lCond}, \eqref{cond:l'>lCond} and \eqref{cond:invertAl}. Then it  is a basis of successive minima for $\Lambda$.
\end{proposition}

\begin{proof}
    First, we prove that for every $1\leq j\leq k$, with $j\in I_{\ell}$, we have $\Vert \mathbf{u}_j\Vert=\mu_{\ell}$. Let $1 \leq \ell \leq m$ and $j \in I_\ell$. Since ${\bf v}_1,\dots,{\bf v}_k$ are orthogonal, by the Lemma \ref{lem:UMIneq}, we have
    $$
    \Vert {\bf u}_j\Vert=\max_{1 \leq \ell' \leq m} \max_{i \in I_{\ell'}} |a_{i,j}|\Vert v_i \Vert.
    $$
    By conditions \eqref{cond:I_ellCond}, \eqref{cond:l'<lCond} and \eqref{cond:l'>lCond}, we have $\Vert \mathbf{u}_j\Vert\leq \mu_{\ell}$. Moreover, by Remark \ref{lem:a_sigma(p),p!=0}, there exists $i_0\in I_{\ell}$ such that $a_{i_0,j}\in \mathbb{F}_q\setminus \{0\}$, so:
    \begin{equation}
        \Vert \mathbf{u}_j\Vert=\vert a_{i_0,j}\vert\times \Vert \mathbf{v}_{i_0}\Vert=\lambda_{i_0}=\mu_{\ell}.
    \end{equation}
    Next, we prove that $\mathbf{u}_1,\dots,\mathbf{u}_k$ are $\mathcal{K}$-linearly independent. Suppose there exist $\{t_j\}_{1 \leq j \leq k} \subset \mathcal{K}$ such that $\sum_{j=1}^kt_j\mathbf{u}_k=0$. Then 
    $$0=\sum_{j=1}^kt_j\mathbf{u}_j=\sum_{\ell=1}^m\sum_{j\in I_{\ell}}t_j\mathbf{u}_j=\sum_{\ell=1}^m\sum_{j\in I_{\ell}}t_j\left(\sum_{i\in I_{\ell}}a_{i,j}\mathbf{v}_i+\sum_{\ell'<\ell}\sum_{ i\in I_{\ell'}}a_{i,j}\mathbf{v}_i\right) \quad \quad \textrm{by \eqref{cond:l'>lCond}}$$
    $$=\sum_{\ell=1}^m\sum_{i\in I_{\ell}}\left(\sum_{j\in \bigcup_{\ell'\geq \ell}I_{\ell'}}a_{i,j}t_j\right)\mathbf{v}_i.$$
    By the linear independence of $\mathbf{v}_1,\dots, \mathbf{v}_k$, we get $\sum_{j\in \bigcup_{\ell'\geq \ell}I_{\ell'}}a_{i,j}t_j=0$, for all $1 \leq \ell \leq m$ and $i\in I_{\ell}$. In particular, for every $i\in I_m$, we have $\sum_{j\in I_m}a_{i,j}t_j=0$. Hence, by (\ref{cond:invertAl}), for every $j\in I_m$, we have $t_j=0$. This implies that
    \begin{equation}
        \sum_{\ell=1}^{m-1}\sum_{i\in I_{\ell}}t_j\mathbf{u}_j=0.    
    \end{equation}
    By a similar argument, we conclude that $t_j=0$ for all $j\in I_{m-1}$. Continuing inductively, we deduce that $t_j=0$, for all $1 \leq  \ell \leq m$ and $j\in I_{\ell}$. Consequently, $\mathbf{u}_1,\dots,\mathbf{u}_k$ are $\mathcal{K}$-linearly independent. By Lemma \ref{rm:succindimplybasss}, we conclude that $({\bf u}_1,\dots,{\bf u}_k)$ is a basis of successive minima for $\Lambda$.
\end{proof}

\begin{proof}[Proof of Theorem \ref{thm:numOrtbases}]
    By Proposition \ref{prop:condnecfororth} and Proposition \ref{prop:condsuffororth}, $(\mathbf{u}_1,\dots,\mathbf{u}_k)$ is an orthogonal basis of $\Lambda$ if and only if the conditions \eqref{cond:I_ellCond}, \eqref{cond:l'<lCond}, \eqref{cond:l'>lCond}, and \eqref{cond:invertAl} are satisfied. In particular, for each  $j\in I_{\ell}$, we can express the vector $\mathbf{u}_j$ as
    $$\mathbf{u}_j=\sum_{\ell'<\ell}\sum_{i\in I_{\ell'}}a_{i,j}\mathbf{v}_i+\sum_{i\in I_{\ell}}a_{i,j}\mathbf{v}_i.$$ 
     Firstly, for every $j\in I_{\ell}$, by \eqref{cond:l'<lCond}, there are $\prod_{\ell'<\ell}\left(q\frac{\mu_{\ell}}{\mu_{\ell'}}\right)^{s_{\ell'}}$ choices for $\{a_{i,j}\}_{i\in \bigcup_{\ell'<\ell}I_{\ell'}}\subseteq \mathcal{R}$. Next, by \eqref{cond:invertAl} and \cite[Equation (2) on page 35]{DummitFoote2004}, there are $|\operatorname{GL}_{s_{\ell}}(\mathbb{F}_q)|=\prod_{i=0}^{s_{\ell}-1}(q^{s_{\ell}}-q^i)$ ways to choose $\{ a_{i,j}\}_{i, j \in I_\ell}$. 
     
     Thus, if $m=1$, then, $s_1=k$, and so there are $\prod_{i=0}^{k-1}(q^k-q^i)$ orthogonal ordered bases for $\Lambda$. If $m\geq 2$, then, there are
    \begin{equation}
        \begin{split}
        \prod_{\ell=1}^m\left(\prod_{\ell'<\ell}\left(q\frac{\mu_{\ell}}{\mu_{\ell'}}\right)^{s_{\ell'}s_{\ell}}\times \prod_{i=0}^{s_{\ell}-1}(q^{s_{\ell}}-q^{i})\right)\\
        =q^{\sum_{\ell=1}^ms_{\ell}\sum_{j=1}^{\ell-1}s_j}\prod_{\ell=1}^m\left(\mu_{\ell}^{s_{\ell}\left(\sum_{i=1}^{\ell-1}s_i-\sum_{i=\ell+1}^ms_i\right)}\prod_{i=0}^{s_{\ell}-1}(q^{s_{\ell}}-q^i)\right)
        \end{split}
    \end{equation}
    orthogonal ordered bases for $\Lambda$. 
\end{proof}

\subsection{Proof of Theorem \ref{thm:NumG_n,k,q(mu,s)}}
The key idea behind counting $(\boldsymbol{\mu},\boldsymbol{s})$-sublattices is to first count the number of orthogonal bases which generate a $(\boldsymbol{\mu},\boldsymbol{s})$-sublattice, and then to use Theorem \ref{thm:numOrtbases} to account for orthogonal bases which generate the same $(\boldsymbol{\mu},\boldsymbol{s})$-sublattice.
\begin{proof}[Proof of Theorem \ref{thm:NumG_n,k,q(mu,s)}]
By Theorem \ref{thm:AraBase}, a $(\boldsymbol{\mu},\boldsymbol{s})$-sublattice is spanned over $\mathcal{R}$ by an orthogonal basis of successive minima. Hence, we first count the number of such bases which generate a lattice $\Lambda\in G_{n,k,q}(\boldsymbol{\mu},\boldsymbol{s},\mathcal{R})$. By Proposition \ref{prop:OrtIndConn}, we can construct an orthogonal set generating a $(\boldsymbol{\mu},\boldsymbol{s})$-sublattice in the following way:
\begin{enumerate}
\item \label{cond:IndInF_q^n} Choose an independent set $\{\mathbf{u}_1,\dots,\mathbf{u}_k\} \subset \mathbb{F}_q^n$.
\item \label{cond:ResOn<} For each $i \in I_{\ell}$, choose $\{{\bf w}_{i} \mid i \in I_\ell \} \subset \mathcal{R}^n_{\leq \log_q(\mu_{\ell})-1}$. 
\item \label{cond:LiftToK^n}Define:
$$\mathbf{v}_i=x^{\log_q(\mu_{\ell})}\mathbf{u}_i+\mathbf{w}_i,$$
for all $1 \leq \ell \leq m$ and $i\in I_{\ell}$. 
\end{enumerate}
By Proposition \ref{prop:OrtIndConn}, \eqref{cond:IndInF_q^n} ensures that tuple $({\bf v}_1,\dots,{\bf v}_k)$ is orthogonal. Moreover, \eqref{cond:ResOn<} and \eqref{cond:LiftToK^n} ensure that $\Lambda=\operatorname{span}_{\mathcal{R}}\left\{\mathbf{v}_1,\dots,\mathbf{v}_k\right\}$ is a $(\boldsymbol{\mu},\boldsymbol{s})$-sublattice. To count the number of such ordered vectors, observe the following:
\vskip 1mm
\noindent $\bullet$ There are $\prod_{i=0}^{k-1}(q^n-q^i)$ ways to choose an independent set $\{\mathbf{u}_1,\dots,\mathbf{u}_k\}\subset \mathbb{F}_q^n$.
\vskip 1mm
\noindent $\bullet$ For each $i\in I_{\ell}$, there are $\mu_{\ell}^n$ ways to choose $\mathbf{w}_i \in \mathcal{R}^n_{\leq \log_q(\mu_{\ell})-1}$. 
\vskip 1mm
\noindent
Hence, the number of ordered $k$-tuples which generate a $(\boldsymbol{\mu},\boldsymbol{s})$-sublattice is:
\begin{equation}
\label{eqn:kTupleCount}
\begin{cases}
    \mu_1^{kn}\prod_{i=1}^{k-1}(q^n-q^i)&\text{if }m=1,\\
    \prod_{i=0}^{k-1}(q^n-q^i)\prod_{\ell=1}^m\mu_{\ell}^{ns_{\ell}}&\text{otherwise}.
\end{cases}
\end{equation}
Hence, to obtain the number of $(\boldsymbol{\mu},\boldsymbol{s})$-sublattices, we divide (\ref{eqn:kTupleCount}) by the number of orthogonal bases for a fixed $(\boldsymbol{\mu},\boldsymbol{s})$-sublattice. Thus, by Theorem \ref{thm:numOrtbases},
\begin{equation}
    \vert G_{n,k,q}(\boldsymbol{\mu},\boldsymbol{s},\mathcal{R})\vert=\begin{cases}
        \mu_1^{kn}&\textrm{if}\,\,m=1,\\
        \frac{\prod_{i=0}^{k-1}(q^n-q^i)}{\prod_{\ell=1}^m\prod_{i=0}^{s_{\ell}-1}(q^{s_{\ell}}-q^i)}\frac{\prod_{\ell=1}^m\mu_{\ell}^{s_\ell\left(n-\sum_{i=1}^{\ell-1}s_i+\sum_{i=\ell+1}^ms_i\right)}}{q^{\sum_{1 \leq j<\ell \leq m}s_js_{\ell}}}&\text{otherwise}.
    \end{cases}
\end{equation}
\end{proof}
\subsection{Proof of Theorem \ref{thm:PrimLattCount}}
In order to count the number of primitive $(\boldsymbol{\mu},\boldsymbol{s})$-sublattices, we first provide an equivalent definition of primitivity.

\begin{definition}
    Let $\mathbf{v}_1,\dots,\mathbf{v}_k\in \mathcal{R}^n \setminus \{{\bf 0}\}$. We say that a monic polynomial $D\in \mathcal{R}$ is \emph{the greatest common divisor} of $\mathbf{v}_1,\dots,\mathbf{v}_k$ if 
    \begin{enumerate}
        \item \label{cond:Ddivisor} for every $1 \leq i \leq k$, there exists $\mathbf{u}_i\in \mathcal{R}^n$, such that $\mathbf{v}_i=D\mathbf{u}_i$;
        \item for every $D'\in \mathcal{R}$ satisfying (\ref{cond:Ddivisor}), we have $D'|D$.
    \end{enumerate}
    If a polynomial $D$ satisfies both conditions, we write $D=\gcd(\mathbf{v}_1,\dots,\mathbf{v}_k)$.  
\end{definition}

\begin{lemma}
\label{lem:primGCD}
    Let $\Lambda \leq \mathcal{R}^n$ be a $k$-sublattice with $\mathcal{R}$-basis $\{\mathbf{v}_1,\dots,\mathbf{v}_k\}$ . Then $\Lambda$ is a primitive if and only if $\operatorname{gcd}(\mathbf{v}_1,\dots,\mathbf{v}_k)=1$. 
\end{lemma}
\begin{proof}
This is straightforward.
\end{proof}
We also need a lemma about sums of the M\"{o}bius function in $\mathcal{K}$. Towards this end,  a monic polynomial $P\in \mathcal{R}$ is called \emph{irreducible} if the only monic polynomials dividing $P$ are $1$ and $P$. Define the M\"{o}bius function $\mu:\mathcal{R}\rightarrow \{-1,0,1\}$ by
$$\mu(N)=\begin{cases}
    (-1)^{\ell}&N=\prod_{i=1}^{\ell}P_i\text{  where }P_i\text{  are distinct irreducible polynomials}\\
    0&\text{else}
\end{cases}.$$
\begin{lemma}{\cite[Lemma 2.1]{BCJ}}
\label{lem:Mobius}
    Let $q$ be a prime power and let $d\in \mathbb{N}$. Then, 
    $$\sum_{N\in \mathcal{R},\deg(N)=d,N\text{ is monic}}\mu(N)=\begin{cases}
        1&d=0\\
        -q&d=1\\
        0&\text{else}
    \end{cases}.$$
\end{lemma}
\begin{proof}[Proof of Theorem \ref{thm:PrimLattCount}]
 By Lemma \ref{lem:primGCD} and the inclusion-exclusion principle, the number of primitive $(\boldsymbol{\mu},\boldsymbol{s})$-sublattices is:
\begin{equation}
\begin{split}
    \vert \widehat{G}_{n,k,q}(\boldsymbol{\mu},\boldsymbol{s},\mathcal{R})\vert=\left \vert \left(\left(G_{n,k,q}(\boldsymbol{\mu},\boldsymbol{s},\mathcal{R})\right)\setminus \bigcup_{P\in \mathcal{R}, P\text{  prime}}P\mathcal{L}_{n,k}(\mathcal{R})\cap G_{n,k,q}(\boldsymbol{\mu},\boldsymbol{s},\mathcal{R})\right) \right \vert\\
    =|G_{n,k,q}(\boldsymbol{\mu},\boldsymbol{s},\mathcal{R})|+\sum_{P_1,\dots,P_{\ell}\text{  prime}}(-1)^{\ell}|P_1\cdots P_{\ell}\mathcal{L}_{n,k}(\mathcal{R})\cap G_{n,k,q}(\boldsymbol{\mu},\boldsymbol{s},\mathcal{R})|\\
    =\sum_{N\in \mathcal{R}, N\text{  is monic}}\mu(N)|N\mathcal{L}_{n,k}(\mathcal{R})\cap G_{n,k,q}(\boldsymbol{\mu},\boldsymbol{s},\mathcal{R})|,
\end{split}
\end{equation}
where $\mu$ is the M\"{o}bius function. Note that if $N|\gcd(\mathbf{v}_1,\dots,\mathbf{v}_k)$, then $\Lambda=\operatorname{span}_{\mathcal{R}}\{\mathbf{v}_1,\dots,\mathbf{v}_k\}$ is a $(\boldsymbol{\mu},\boldsymbol{s})$-sublattice, if and only if $\frac{1}{N}\Lambda$ is a $\left(\frac{\boldsymbol{\mu}}{\vert N\vert},\boldsymbol{s}\right)$-sublattice in $\mathcal{R}^n$. Therefore, by Lemma \ref{lem:Mobius}, we have
\begin{equation}
\begin{split}
\label{eqn:GhatComp}
    \vert\widehat{G}_{n,k,q}(\boldsymbol{\mu},\boldsymbol{s},\mathcal{R})\vert&=\sum_{N\in \mathcal{R},N\text{  is monic}}\mu(N)\left \vert G_{n,k,q}\left(\frac{\boldsymbol{\mu}}{\vert N\vert},\boldsymbol{s},\mathcal{R}\right) \right \vert\\
    &=\sum_{d=0}^{\infty}\left \vert G_{n,k,q}\left(\frac{\boldsymbol{\mu}}{q^d},\boldsymbol{s},\mathcal{R}\right) \right \vert \times \left(\sum_{N\in \mathcal{R},\deg(N)=d,N\text{  is monic}}\mu(N)\right)\\
    &=|G_{n,k,q}\left(\boldsymbol{\mu},\boldsymbol{s},\mathcal{R}\right)|-q\left\vert G_{n,k,q}\left(\frac{\boldsymbol{\mu}}{q},\boldsymbol{s},\mathcal{R}\right)\right\vert.
\end{split}
\end{equation}

Suppose first that $\mu_1=1$. Then, $G_{n,k,q}\left(\frac{\boldsymbol{\mu}}{q},\boldsymbol{s},\mathcal{R}\right)=\emptyset$. Hence, by (\ref{eqn:GhatComp}), $\widehat{G}_{n,k,q}(\boldsymbol{\mu},\boldsymbol{s},\mathcal{R})=G_{n,k,q}(\boldsymbol{\mu},\boldsymbol{s},\mathcal{R})$. 

Suppose now that $m=1$ and $\mu_1\neq 1$. Then, by Theorem \ref{thm:NumG_n,k,q(mu,s)} and (\ref{eqn:GhatComp}), $\vert \widehat{G}_{n,k,q}(\boldsymbol{\mu},\boldsymbol{s},\mathcal{R})\vert=\mu_1^{nk}\left(1-q^{-(nk-1)}\right)$. 

Suppose finally that $\mu_1 \neq 1$ and $m \geq 2$. Then, by Theorem \ref{thm:NumG_n,k,q(mu,s)} and (\ref{eqn:GhatComp}), we obtain: 
$$\vert \widehat{G}_{n,k,q}(\boldsymbol{\mu},\boldsymbol{s},\mathcal{R})\vert=\frac{\prod_{i=0}^{k-1}(q^k-q^i)\times \prod_{j=1}^m\left(\mu_j^{s_j\left(n-\sum_{i=1}^{\ell-1}s_i+\sum_{i=\ell+1}^ms_i\right)}\times \left(1-q^{1-s_j\left(n-\sum_{i=1}^{\ell-1}s_i+\sum_{i=\ell+1}^ms_i\right)}\right)\right)}{q^{\sum_{\ell<j}s_js_{\ell}}\prod_{\ell=1}^m\prod_{i=0}^{s_{\ell}-1}(q^{s_{\ell}}-q^i)}.$$
\end{proof}
\bibliography{Ref}
\bibliographystyle{amsalpha}
\end{document}